\documentclass[11pt,reqno]{amsart}
\usepackage[colorlinks=true, linkcolor=blue]{hyperref}
\numberwithin{equation}{section}
\usepackage{times}
\usepackage{amsmath,amsfonts,amstext,amssymb,amsbsy,amsopn,amsthm,eucal}
\usepackage{txfonts}
\usepackage{dsfont}
\usepackage{graphicx}   % for figures
\usepackage{accents}
\usepackage{enumerate}
\usepackage{xcolor}
\usepackage{verbatim}
\usepackage{esint}

\usepackage[normalem]{ulem}
\usepackage{cancel}

\newcommand{\R}{\text{R}}

\newtheorem{theorem}{Theorem}[section]

\newtheorem{proposition}[theorem]{Proposition}
\newtheorem{lemma}[theorem]{Lemma}
\newtheorem{corollary}[theorem]{Corollary}
\theoremstyle{definition}
\newtheorem{definition}[theorem]{Definition}

\theoremstyle{remark}
\newtheorem{remark}[theorem]{Remark}
\theoremstyle{remark}

\theoremstyle{remark}

\theoremstyle{remark}
\theoremstyle{remark}

\title{Singular sets of homogeneous Dirichlet heat flows at Large time}

\author{Xiujin Chen$^{*}$ and Junyuan Wang$^{*}$}

\thanks{School of Mathematical Sciences, Zhejiang University,
Hangzhou 310058, China.}

\thanks{$^{*}$Corresponding authors:
Xiujin Chen (\texttt{xjchen1998@foxmail.com}) and
Junyuan Wang (\texttt{wangjunyuan@zju.edu.cn}).}

\begin{document}

\maketitle

\vspace{-10mm}

\begin{abstract}
    We study homogeneous Dirichlet heat flows for real symmetric, uniformly elliptic divergence-form operators. We establish quantitative stability of singular sets near homogeneous harmonic polynomials and eigenfunctions. Based on these, in the smooth coefficient case, we prove the spatial singular set at sufficiently large time is codimension-two and has a uniform measure upper bound independent of time. However, this codimension-two property fails in the Lipschitz coefficient case. On a three-dimensional ball, we construct a fixed self-adjoint operator and a heat flow consisting of two eigenmodes whose singular sets contain two-dimensional disks at a sequence of times tending to infinity. 
\end{abstract}

\vspace{8pt}

% --------------------------------------------------------------------
% Problem, definitions, and main results
\section{Introduction}\label{problem-definitions-and-main-results}
Let \(n\ge2\), and \(\Omega\subset\mathbb R^n\) be a bounded
connected domain with \(\partial\Omega\in C^\infty\). In this paper, we consider the following homogeneous Dirichlet heat flow,
\begin{equation}\label{Dirichlet problem}
\begin{cases}
\partial_tu=\operatorname{div}(A(x)\nabla u)
 &\text{in }\Omega\times(0,\infty),\\
u=0&\text{in }\partial\Omega\times(0,\infty).
\end{cases}
\end{equation}
The matrix \(A\) is real symmetric and independent of time. Assume that, for some \(\kappa>0\),
\begin{equation}\label{uniformly elliptic}
\kappa|\xi|^2\le\xi^TA(x)\xi\le\kappa^{-1}|\xi|^2.
\end{equation}
\indent Throughout this paper, solutions are real-valued Dirichlet heat flows in the standard energy class. We assume nonzero initial data \(u(0)=u_{0}\in L^2(\Omega)\) and that the solution subsequently
evolves under the self-adjoint heat semigroup. Define the interior singular set of a time slice by
\begin{equation}
S_t:=S_{t}(u)=\{x\in\Omega:u(x,t)=0,\ \nabla_xu(x,t)=0\},
\end{equation}
and the global singular set by
\begin{equation}
 \Sigma_t=\{x\in\overline{\Omega}:u(x,t)=0,\ \nabla_xu(x,t)=0\}.
\end{equation}
\noindent This definition involves only the spatial gradient; it does not require \(\partial_tu=0\). Let \(L=-\operatorname{div}(A\nabla)\) and \(P_\lambda\) be the spectral projection onto the corresponding eigenspace. Set
\begin{equation*}
\lambda_*=\min\{\lambda:P_\lambda u_{0}\ne0\},\qquad
\Phi=P_{\lambda_*}u_{0}\ne0.
\end{equation*}
%We call \(\Phi\) the leading eigenfunction and set
%\begin{equation*}
%K_\Phi=\{p\in\overline\Omega:\Phi(p)=0,\ \nabla\Phi(p)=0\}.
%\end{equation*}
\indent An important fact is that scalar normalization selects a nonzero elliptic eigenfunction as its large-time limit, i.e.,
\begin{equation*}
    e^{\lambda_{*}t}u(\cdot,t)\to \Phi\;\text{as}\;t\to\infty.
\end{equation*}
%Vanishing orders here include boundary points. When sufficient
%regularity is available, the vanishing order is the degree of the first
%nonzero homogeneous term in the Taylor expansion.
Classical result tells us the singular set of $\Phi$ is codimension two. Hence, a natural question is whether the spatial singular sets inherit the codimension-two geometry of this limit at every sufficiently large time, i.e., \textit{for each fixed nonzero solution, does there exist \(T_0\) such that \(\dim_HS_t\le n-2\) for every \(t>T_0\)?} 

 Actually, this question is raised in \cite{JW1}. In \cite{JW1}, using blow-up analysis, the authors gave a new proof of the main result of one-dimensional parabolic equations in \cite{Angenent}, as a corollary, the spatial singular sets are empty at large time. This inspired the posing of the above higher-dimensional question.

In the following, we first review the main developments in the study of singular sets of elliptic equations. In 1994, Han \cite{Hansingular} proved the $(n-2)$-rectifiablity for singular sets of elliptic equations with H\"older coefficients. In \cite{98singular}, Han, Hardt and Lin proved the $\mathcal{H}^{n-2}$-finiteness of singular sets of the solutions to elliptic equations with smooth coefficients. Later, the $\mathcal{H}^{n-2}$-finiteness of critical sets is proved by Hardt, Hoffmann-Ostenhof, M, Hoffmann-Ostenhof, T and Nadirashvili \cite{hardt1999} under the same setting. For non-smooth coefficients, Cheeger, Naber and Valtorta \cite{CNV} introduced the quantitative stratification technique to show the weak volume estimates for critical sets under Lipschitz coefficients. Later, Naber and Valtorta \cite{NV} improved this result to a uniform Minkowski-type volume estimate. For more weak coefficients, Kim  \cite{Kim} generalized the structure theorem in \cite{Hansingular} to quasilinear elliptic equations with Dini coefficients. In 2022, Kenig and Zhao \cite{KZ1} obtained the volume estimates for singular sets of harmonic functions in $C^{1}$-Dini domain. More recently, the results in \cite{NV} were generalized to H\"older case by Huang and Jiang \cite{HJ1} and Dini case by Jiang and Wang \cite{JW2}  under the doubling assumption.  

 Concerning the geometric sets at time slices of solutions to parabolic equations, since $u(\cdot,t)$ does not satisfy the usual elliptic equation, hence it brings some difficulties to analyze its geometric sets. In 1991, Lin \cite{Linconj} considered the heat equation in an analytic closed Riemannian manifold, and he showed a sharp upper bound for spatial nodal set. Recently, Huang and Jiang \cite{HJ2} proved the $\mathcal{H}^{n-1}$-finiteness of nodal set at time slice for Lipschitz coefficients. Also in Lipshcitz case, Hallgrent-Koirala-Ma \cite{HKM} proved a uniform parabolic Minkowski estimate for interior singular set at space-time, and as a corollary, they claimed that the spatial singular sets have a uniform $\mathcal{H}^{n-2}$ measure bound for $\mathcal{H}^{1}$-almost every time.
 
  In our case, although the limiting eigenfunction has a codimension-two singular set, the convergence of the functions alone does not imply that the spatial singular sets of solutions share the same geometric property. Therefore, we need to prove some stability result in order to obtain a comparison of the critical sets. In fact, if the coefficient matrix $A(x)$ is sufficiently smooth, we can give a positive answer to the initial question. The following is the first main theorem in this paper.
\begin{theorem}
\label{thm:pos-main}
Assume that $A(x)\in C^{\infty}(\overline{\Omega})$ and let $u$ be a nonzero solution to \eqref{Dirichlet problem} with \eqref{uniformly elliptic}. Under the above setting, there exist positive constants $T=T(n,A,\Omega,u_{0})$ and 
 $C=C(n,A,\lambda_{*},\Omega)$
such that
\begin{equation}\label{eq:pos-main-measure}
 \sup_{t\ge T}\mathcal{H}^{n-2}(S_t)
 \le\sup_{t\ge T}\mathcal{H}^{n-2}(\Sigma_t)\le C<\infty.
\end{equation}
Especially, for each $t\geq T$, we have $\dim_{\mathrm H}\Sigma_t\le n-2$.
\end{theorem}
\begin{remark}
    Here, we proved the desired result for global singular set $\Sigma_{t}$, which includes the boundary singularities, and the measure bound is independent of time.
\end{remark}
\indent The proof has three parts. Spectral normalization makes the heat
flow converge in every finite spatial \(C^M\) norm to a fixed
Dirichlet eigenfunction. A quantitative perturbation lemma gives
a uniform measure bound near each finite-order harmonic model.
Finally, boundary finite order and a finite covering turn the
local estimates into a bound on the closed domain.
Compactness of the unit sphere in each eigenspace makes the
bound independent of the particular limiting eigenfunction.

 Analogous to the singular sets of elliptic equations, it is natural to ask whether the above theorem holds for Lipschitz case. Unfortunately, for coefficients with low regularity, we have the following counterexample, which is the second main theorem.

\begin{theorem}\label{thm:counter}
On \(\Omega=B_1(0)\subset\mathbb{R}^{3}\), there exist a real symmetric, time-independent, uniformly elliptic Lipschitz matrix \(A\) and a nonzero heat flow \(u\) satisfying \eqref{Dirichlet problem} and \eqref{uniformly elliptic}, such that along a sequence \(t_j\to\infty\), the spatial singular set \(S_{t_j}(u)\) contains the diffeomorphic image of a two-dimensional disk. Therefore,
\begin{equation*}
\dim_HS_{t_j}(u)\ge2>1=n-2.
\end{equation*}
\end{theorem}
\begin{remark}
    In the construction, $u$ consists of just two Dirichlet eigenmodes, and the final matrix $A(x)$ is permitted to belong to any prescribed class \(C^{1,\alpha}\) with \(0<\alpha<1/2\).
\end{remark}
\begin{remark}
    Actually, for time-dependent coefficient matrix $A(x,t)$, although we assume that it is $C^{\infty}$-smooth in space-time, one can also construct a corresponding counterexample in a similar way in $\mathbb{R}^{4}$.
\end{remark}
\indent The distinction concerns the stability of singular sets near the limiting spectral profile. In the smooth setting, convergence in sufficiently high spatial norms places the entire normalized tail inside a fixed geometric stability neighborhood. The counterexample preserves low-order regularity while allowing higher derivatives to deteriorate on shrinking spatial scales, producing critical level sets that are selected by the decaying second mode.

\noindent \textbf{Disclosure of AI tools.} 
For the smooth case, the main idea of the proof originated from the authors. By using scalar normalization to converge to an eigenfunction, together with the method for estimating singular sets of smooth elliptic coefficients, we came up with the idea of proving Lemma 3.9, which was already known for two-dimensional harmonic functions. We therefore attempted to use GPT 6.0 Sol to generalize it to higher dimensions. For the Lipschitz coefficient case, after some attempts, we realized that obtaining a positive result would be rather difficult, so we turned to GPT 6.0 Sol to construct a counterexample. However, the first draft was poorly readable. Through multiple rounds of interactive discussion and supplementation, we reorganized and rewrote the propositions, arriving at the current version. The construction of the counterexample now has a clearer logic.

\noindent \textbf{Acknowledgements.}
  This problem was raised by Jiang and the second author in \cite{JW1}. The authors thank Prof. Wenshuai Jiang for reminding us of it again and informed us that this should be true in the analytic setting. We are also grateful to him for his constant help and encouragement. This work was supported by the National Key Research and Development Program of China (No. 2022YFA1005500).

\part*{\texorpdfstring{\centering\Large Part I. Codimension-Two Spatial Singular Set in the Smooth Case}{Part I. Codimension-Two Spatial Singular Set in the Smooth Case}}
\indent In the first part, we assume that $A(x)$ is a smooth coefficient matrix. Part I is organized as follows. We recall some needed basic properties in Section 2. Section 3 prepares the required good two-planes, uniform fiber cardinality and then we show the critical set comparison lemma. In Section 4, we deal with the boundary and prove Theorem \ref{thm:pos-main}.
\section{Preliminaries}
\begin{lemma}\cite[Sections 4.8-4.12]{pos-Hunter}.\label{lem:pos-spectrum}
The quadratic form
\[
 \mathfrak a(v,w)=
 \int_{\Omega}\sum_{i,j=1}^n
 a_{ij}\partial_{x_j}v\,\partial_{x_i}w\,dx,
 \qquad v,w\in H_0^1(\Omega),
\]
defines a positive self-adjoint operator \(L\) on \(L^2(\Omega)\).
Its resolvent is compact. Its distinct eigenvalues can be listed as
\[
 0<\lambda_1<\lambda_2<\cdots,\qquad \lambda_j\to\infty,
\]
with finite-dimensional mutually orthogonal eigenspaces spanning
\(L^2(\Omega)\). Eigenfunctions are smooth on \(\overline{
\Omega}\).
The heat flow is \(u(t)=e^{-tL}u_0\), with
\(L=-{\rm div}(A\nabla)\) and homogeneous Dirichlet boundary condition.
\end{lemma}

\begin{lemma}\label{lem:pos-smoothing}
For every \(s>0\) and integer \(M\ge0\),
\[
 e^{-sL}:L^2(\Omega)\longrightarrow C^M(\overline{\Omega})
\]
is bounded.
\end{lemma}

\begin{lemma}[The lowest occupied spectral profile]\label{lem:pos-profile}
Let \(P_\lambda\) denote orthogonal projection onto the entire
eigenspace for a distinct eigenvalue \(\lambda\).
Set
\[
 \lambda_*=\min\{\lambda:P_\lambda u_0\ne0\},\qquad
 \Phi=P_{\lambda_*}u_0.
\]
Then \(\Phi\ne0\), \(L\Phi=\lambda_*\Phi\), and
\(\Phi|_{\partial\Omega}=0\).
If
\[
 w(t)=e^{\lambda_*t}u(t)=\Phi+R(t),
\]
then there is \(\delta>0\) such that for every integer \(M\ge0\),
\begin{equation}\label{eq:pos-high-convergence}
 ||{R(t)}||_{C^M(\overline{\Omega})}\le C_M e^{-\delta t},
 \qquad t\ge1.
\end{equation}
\end{lemma}

\begin{proof}
At least one spectral projection of \(u_0\) is nonzero, so the
minimum is well-defined. Retaining the entire eigenspace handles
eigenvalue multiplicity without choosing a preferred eigenvector.
The spectral expansion is
\[
 R(t)=\sum_{\lambda>\lambda_*}
 e^{-(\lambda-\lambda_*)t}P_\lambda u_0.
\]
Let \(\delta\) be the gap from \(\lambda_*\) to the next distinct
eigenvalue. Orthogonality yields
\[
 ||{R(t)}||_{L^2}^2
 =\sum_{\lambda>\lambda_*}
 e^{-2(\lambda-\lambda_*)t}||P_\lambda u_0||_{L^2}^2
 \le e^{-2\delta t}||{u_0}||_{L^2}^2.
\]
For \(t\ge1\),
\[
 R(t)=e^{\lambda_*/2}e^{-L/2}R(t-1/2).
\]
Lemma~\ref{lem:pos-smoothing} now gives
\eqref{eq:pos-high-convergence}, after absorbing fixed factors
into \(C_M\).
\end{proof}
\indent Next, we introduce some algebraic background that will be needed. First is the concept of Hilbert series.
\begin{definition}
    For graded space $M=\bigoplus_{k\geq 0}M_{k}$, its Hilbert series is defined by 
    \begin{equation*}
        H_{M}(t)=\sum_{k\geq 0}({\rm dim_{\mathbb{C}}}M_{k})t^{k}.
    \end{equation*}
\end{definition}
Let $R=\mathbb{C}[y_{1},y_{2}]$ be the complex polynomial ring, $f$ and $g$ are relatively prime homogeneous polynomials  of degree $d$ with real coefficients.
\begin{lemma}\label{lem2.5}
    Let $A=R/(f,g)$ be the quotient algebra, then ${\rm dim}_{\mathbb{C}}A=d^{2}$.
\end{lemma}
\begin{proof}
    \indent Denote $R_{k}$ be the homogeneous polynomial space of degree $k$, then $R=\bigoplus_{k}R_{k}$. Set $R_{k}=0$ for any $k\leq 0$, and for any $s\in\mathbb{N}$, let $R(-s)_{k}:=R_{k-s}$. Note that $f$ and $g$ are relatively prime, then every relation \(af+bg=0\) consequently has the form
\(a=gc\), \(b=-fc\). We define the following maps.
\begin{align*}
\begin{split}
    \alpha:R\rightarrow R\oplus R&,\quad\alpha(c)=(gc,-fc),\\
     \beta: R\oplus R \rightarrow R&,\quad\beta(a,b)=fa+gb,\\
     q: R\rightarrow A&,\quad q(h)=[h].
    \end{split}
\end{align*}Therefore the graded sequence
\[
 0\rightarrow R(-2d)
 \xrightarrow{\alpha} R(-d)\oplus R(-d)
 \xrightarrow{\beta} R
 \xrightarrow{q} A\longrightarrow0
\]
is exact. By standard calculation, we can see its Hilbert series is
\[
 \frac{1-2t^d+t^{2d}}{(1-t)^2}
 =(1+t+\cdots+t^{d-1})^2.
\]
Hence the complex dimension is \(d^2\), and all homogeneous
components of degree at least \(N_0=2d-1\) vanish. The generators are real, so the real polynomial quotient has real dimension \(d^2\).
\end{proof}
\section{Quantitative stability of critical sets near homogeneous harmonic polynomials}
\label{sec:pos-geometry}
\subsection{Simultaneously good two-planes and finite local algebras}

For \(1\le i<j\le n\), let
\[
 \pi_{ij}:\mathbb{R}^n\longrightarrow\mathbb{R}^{n-2}
\]
delete coordinates \(i\) and \(j\), and let
\(\iota_{ij}:\mathbb{R}^2\to\mathbb{R}^n\) insert its two arguments into
these coordinates, setting the remaining coordinates to zero.
When \(n=2\), the target of \(\pi_{12}\) is the one-point space.

\begin{lemma}\cite[Corollary 3.2]{hardt1999}
\label{lem:pos-plane}
Let \(P\) be a nonzero real homogeneous harmonic polynomial
of degree \(m\ge2\) on \(\mathbb{R}^n\).
There is one real orthogonal coordinate system in which,
for  \(1\leq i<j\leq n\), $p_{ij}(y):=P(\iota_{ij}y)$ 
has an isolated critical point at \(0\) after complexification
to \(\mathbb{C}^2\).
\end{lemma}

\begin{remark}
For \(n>2\), \(p_{ij}\) need not be harmonic in its two
variables. Its isolated complex critical point is the property
used in the local algebra.
\end{remark}

\begin{lemma}[Finite differentiability preparation theorem]
\label{lem:pos-preparation}
Let \(F_0:\mathbb{R}^N\to\mathbb{R}^{N}\) be smooth with $F(0)=0$. Let
\(\mathcal E_N\) be the ring of smooth real germs at \(0\), and assume
\[
 \mu=\dim_{\R}\mathcal E_N/(F_{0,1},\ldots,F_{0,N})<\infty.
\]
Choose smooth representatives \(e_1,\ldots,e_\mu\) of a basis. Fix a small bounded neighborhood $0\in D\subset \mathbb{R}^{N}$. There are fixed source neighborhoods \(0\in V\subset U\subset D\), a target neighborhood \(0\in Q\), and $\eta>0$ such that for any $G\in C^{2\mu}(D,\mathbb{R}^{N})$ with $||G-F_{0}||_{C^{2\mu}(D)}\leq \eta$, then $G(V)\subset Q$ and each \(a\in C^{2\mu}(U)\) admits
\begin{equation}\label{eq:pos-preparation}
 a(x)=\sum_{l=1}^{\mu}e_{l}(x)a_{l}(G(x)),
 \quad x\in V,\quad a_{l}\in C^1(Q).
\end{equation}
\end{lemma}

\begin{remark}
This is the finite generating representation in
\cite[Lemma 2.3]{pos-HHL}, taking both regularity indices there
equal to \(1\); its original source is \cite{pos-Bakhtin}.
Only \eqref{eq:pos-preparation} and fixed neighborhoods are used.
No open-image property of a real finite map is needed.
\end{remark}

\begin{lemma}
\label{lem:pos-algebra}
In the coordinates of Lemma~\ref{lem:pos-plane}, set
\[
 F_{ij}(x)=(\partial_iP(x),\partial_jP(x),\pi_{ij}x).
\]
For every \(i<j\), its smooth local algebra has dimension
\((m-1)^2\). A basis may be represented by real polynomials
in \(x_i,x_j\), with its first representative equal to \(1\).
\end{lemma}

\begin{proof}
We only need to consider the pair $(1,2)$ and take \(p_{12}\) as in Lemma \ref{lem:pos-plane}, we denote
\(p_1=\partial_{y_1}p_{12}\), \(p_2=\partial_{y_2}p_{12}\). By the choice of coordinate in Lemma \ref{lem:pos-plane}, \(p_1,p_2\) have no common nonzero complex zero, and they are relatively prime in the ring \(\mathbb{C}[y_1,y_2]\). Let $A=\mathbb{C}[y_{1},y_{2}]/(p_{1},p_{2})$, by Lemma \ref{lem2.5}, we have ${\rm dim}_{\mathbb{C}}A=(m-1)^{2}$. Then
\begin{equation*}
    (\frac{\mathbb{R}[y_{1},y_{2}]}{(p_{1},p_{2})})\otimes_{\mathbb{R}}\mathbb{C}\cong \frac{\mathbb{C}[y_{1},y_{2}]}{(p_{1},p_{2})},
\end{equation*}
hence 
\begin{equation*}
    {\rm dim}_{\mathbb{R}}\frac{\mathbb{R}[y_{1},y_{2}]}{(p_{1},p_{2})}=(m-1)^{2}.
\end{equation*}
\indent Actually, this also gives the smooth local
algebra. Taylor's formula for a smooth germ gives
\[
 a(y)=T_{N_0-1}a(y)+
       \sum_{|\alpha|=N_0}y^\alpha a_\alpha(y),
\]
with smooth \(a_\alpha\).
Every degree-\(N_0\) monomial belongs to \((p_1,p_2)\).
Thus every smooth class has a polynomial representative of degree
less than \(N_0\). If such a representative is a smooth combination
of \(p_1,p_2\), truncating the coefficients through degree
\(N_0-m\) gives the same polynomial identity; when this index
is negative, the identity forces the representative to be zero.
The smooth and polynomial quotients are therefore isomorphic.

Write \(x=(y,z)\) in the chosen pair coordinates.
Quotienting by the coordinates of \(z\) amounts to setting \(z=0\),
because
\[
 a(y,z)-a(y,0)=\sum_{\alpha=1}^{n-2}z_\alpha
             \int_0^1\partial_{z_\alpha}a(y,tz)\,dt.
\]
The local algebra of \(F_{ij}\) is therefore precisely
\(\mathcal E_2/(p_1,p_2)\), of dimension \(d^2\).
The nonzero class of \(1\) can be included in a real polynomial basis.
\end{proof}

\subsection{Uniform fiber cardinality and a graph covering}

The next consequence of preparation is the quantitative point
that will enter the measure estimate. It concerns the whole
fiber, not only the separate coordinate values of its points.

\begin{lemma}
\label{lem:pos-fibers}
Suppose \eqref{eq:pos-preparation} holds on \(V\) for a fixed
map \(G\), with \(\mu\) generators, for every polynomial \(a\).
Then
\begin{equation}\label{eq:pos-fiber-cardinality}
 \#\{x\in V:G(x)=\eta\}\le\mu
 \qquad\text{for every }\eta\in Q.
\end{equation}
\end{lemma}

\begin{proof}
If a fiber contained \(\mu+1\) distinct points
\(x^{(1)},\ldots,x^{(\mu+1)}\), define, for each \(r\),
\[
 L_r(x)=\prod_{\substack{1\le s\le\mu+1\\s\ne r}}
 \frac{(x-x^{(s)})\cdot(x^{(r)}-x^{(s)})}
 {|x^{(r)}-x^{(s)}|^2}.
\]
These polynomials satisfy \(L_r(x^{(s)})=\delta_{rs}\).
Apply preparation to each \(L_r\), obtaining coefficients
\(a_{r\ell}\). At the selected fiber points,
\[
 \delta_{rs}=
 \sum_{\ell=1}^{\mu}e_\ell(x^{(s)})a_{r\ell}(\eta).
\]
In matrix notation, we have $I_{\mu+1}= EA$, by rank comparison, this is a contradiction. Hence each fiber with more than \(\mu\) points contains such a finite selection, so this also rules out infinite fibers.
\end{proof}

\begin{lemma}
\label{lem:pos-roots}
Let \(W\subset\mathbb{R}^{n-2}\) be a open domain and denote
\[
 q_i(x,z)=x^d+\sum_{k=0}^{d-1}a_{ik}(z)x^k,
 \qquad a_{ik}\in C^1(W),\quad i=1,2.
\]
Then any subset of
\[
 \{(y_1,y_2,z):z\in W,\ q_1(y_1,z)=q_2(y_2,z)=0\}
\]
is contained in a countable union of \(C^1\) embedded
\((n-2)\)-dimensional submanifolds.
\end{lemma}

\begin{proof}
Partition the set according to the exact multiplicities
\(s_i\in\{1,\ldots,d\}\) of the two roots.
On a fixed multiplicity set, put
\[
 H_i(y,z)=\partial_x^{s_i-1}q_i(y_i,z),\qquad i=1,2.
\]
At each of its points, \(H_1=H_2=0\), and the Jacobian in
\((y_1,y_2)\) is diagonal, with nonzero diagonal entries
\(\partial_x^{s_i}q_i(y_i,z)\).
The \(C^1\) implicit function theorem therefore gives a local
graph \(y=h(z)\) containing the nearby part of that multiplicity
set. Multiplicity need not be constant in an ambient neighborhood.
Second countability gives countably many such neighborhoods for
each multiplicity set. There are only finitely many multiplicity
pairs. For \(n=2\), the parameter space is a point and the
graphs are zero-dimensional.
\end{proof}

Before proving the critical set comparison lemma, we introduce the following integral geometric formula.
\begin{lemma}[\cite{Fe} 3.2.22]
\label{lem:pos-projection}
Let \(E\subset B_r\subset\mathbb{R}^n\) be a Borel set covered by countably many \(C^1\) embedded
\(n-2\)-dimensional submanifolds. Assume that for any \(1\leq i<j\leq n\)
and \(z\in\mathbb{R}^{n-2}\), $\#(E\cap\pi_{ij}^{-1}\{z\})\le N$, then
\begin{equation}\label{eq:pos-projection-bound}
 \mathcal{H}^{n-2}(E)\leq C(n)N r^{n-2}.
\end{equation}
\end{lemma}

\subsection{Critical set comparison lemma}
\indent In \cite{98singular}, Han, Hardt and Lin proved that for homogeneous harmonic polynomial $P$ of two variables and degree $m$, for any $v\in C^{2m^{2}}(B_{1})$, if $|v-P|_{C^{2m^{2}}(B_{1})}$ is small enough, then we have the local $\mathcal{H}^{n-2}$-measure estimate for the critical set of $v$, which can be viewed as a critical set comparison result. For convenience in proving the stability result in the next section, we improve this result to general homogeneous harmonic polynomial in $\mathbb{R}^{n}$ as follows.
\begin{lemma}
\label{lem:pos-harmonic-stability}
Let \(P\) be a nonzero real homogeneous harmonic polynomial
of degree \(m\ge2\) on $\mathbb{R}^{n}$. Let $M=2(m-1)^{2}+1$. There exist \(\varepsilon_P>0\) and \(r_P\in(0,1)\) such that each \(v\in C^M(B_1)\) with
\[
 ||{v-P}||_{C^M(B_1)}<\varepsilon_P
\]
has its critical set in \(B_{r_P}\) covered by countably many \(C^1\) embedded \(k\)-dimensional submanifolds. Moreover,
for any \(0<r\le r_P\),
\begin{equation}\label{eq:pos-harmonic-measure}
 \mathcal{H}^{n-2}(\{\nabla v=0\}\cap B_r)
 \leq C(n) (m-1)^2 r^{n-2}.
\end{equation}
\end{lemma}

\begin{proof}
\textbf{Step 1:  uniformly bound intersections with two-plane slices.} 
For $i<j$ and each pair $(i,j)$, applying the single orthogonal rotation of Lemma~\ref{lem:pos-plane}. Then we define
\[
 G_{ij}(x)=(\partial_i v(x),\partial_jv(x),\pi_{ij}x).
\]
Let $\mu=(m-1)^2$, Lemma~\ref{lem:pos-algebra} gives multiplicity \(\mu\) for \(F_{ij}\), and by assumption
\[
 ||{G_{ij}-F_{ij}}||_{C^{2\mu}(B_1)}
 \le C||{v-P}||_{C^{2\mu+1}(B_1)}.
\]
Then Lemma \ref{lem:pos-preparation} supplies fixed neighborhoods \(V_{ij},U_{ij},Q_{ij}\) and \(\mu\) polynomial generators for each pair. There are finitely many pairs. Choose \(\varepsilon_P\) sufficiently small to make all preparations valid, and one radius \(r_P\) with
\[
 B_{r_P}\subset\bigcap_{i<j}V_{ij}.
\]
The orthogonal change of coordinates only changes the original
\(C^M\)-closeness by a fixed factor. Let \(E=\{\nabla v=0\}\cap B_{r_P}\). If \(x\in E\) and \(\pi_{ij}x=z\), then
\(G_{ij}(x)=(0,0,z)\).
Every such point lies in the same fiber of \(G_{ij}\) in
\(V_{ij}\). By Lemma~\ref{lem:pos-fibers}, we have 
\begin{equation}\label{eq:pos-all-slice-counts}
 \#(E\cap\pi_{ij}^{-1}\{z\})\le\mu
 \qquad\text{for every pair and every }z.
\end{equation}
If the intersection is nonempty, its image lies in \(Q_{ij}\)
because \(G_{ij}(V_{ij})\subset Q_{ij}\).
Thus the fiber lemma applies to all relevant targets.

\textbf{Step 2: prove the graph covering needed for the area formula.}
It suffices here to use the pair \((1,2)\).
Write \(x=(y_1,y_2,z)\), \(G=G_{12}\), and choose generators
\(e_1=1,e_2,\ldots,e_\mu\).
Apply Lemma \ref{lem:pos-preparation} to each fixed polynomial \(y_i e_j\), obtaining
\[
 y_i e_j(x)=\sum_{\ell=1}^{\mu}b_{ij,l}(G(x))e_\ell(x),
 \qquad b_{ij,l}\in C^1(Q_{12}),\quad i=1,2.
\]
Let \(B_i(\eta)\) be the ordinary square matrix with entry
\(b_{ij,l}(\eta)\) in row \(j\), column \(\ell\).
The column vector \(e(x)=(e_1(x),\ldots,e_\mu(x))^T\)
is nonzero because its first entry is \(1\), and
\begin{equation*}
 B_i(G(x))e(x)=y_i e(x).
\end{equation*}
Hence \(\det(y_iI_\mu-B_i(G(x)))=0\).
On the set \(Z=\{x\in V_{12}:\partial_1v=\partial_2v=0\}\),
one has \(G(x)=(0,0,z)\). Set
\begin{equation}
 W=\{z:(0,0,z)\in Q_{12}\},\qquad
 q_i(Y,z)=\det(YI_\mu-B_i(0,0,z)).
\end{equation}
These polynomials are monic of degree \(\mu\) in \(Y\), and have \(C^1\) coefficients in \(z\), and vanish at \((y_i,z)\) for each point of \(Z\).
Lemma~\ref{lem:pos-roots} implies a countable \(C^1\) graph
covering of \(Z\), and therefore of \(E\).

 If \(n\ge3\), consider the set \(E\cap B_r\), it is a Borel set since $\nabla v$ is continuous. Now applying Lemma~\ref{lem:pos-projection} to it with \eqref{eq:pos-all-slice-counts}, we can obtain \eqref{eq:pos-harmonic-measure}.
For \(n=2\), the only projection has a one-point target,
so Step 2 directly gives \(\#E\le\mu\), the same bound
with \(r^0=1\). Finally rotate back, orthogonal maps preserve balls,
Hausdorff measure and embedded submanifolds.
\end{proof}

\section{Global stability up to the boundary and the estimate at large time}
\label{sec:pos-global}

We now fix the nonzero Dirichlet eigenfunction \(\Phi\) from
Lemma~\ref{lem:pos-profile}, and write its equation as
\[
 {\rm div}(A\nabla\Phi)+\lambda_*\Phi=0.
\]
\subsection{Finite vanishing order and boundary extension}
In this subsection, for the completeness, we show the vanishing order of eigenfunction at each point in $\overline{\Omega}$ is finite.
\begin{lemma}
\label{lem:pos-finite-order}
Let
\[
 K_\Phi=\{p\in\overline{\Omega}:\Phi(p)=0,\ \nabla\Phi(p)=0\}.
\]
At every \(p\in K_\Phi\), the smooth one-sided Taylor expansion,
in local coordinates centered at \(p\), has a first nonzero
homogeneous term \(P_p\) of finite degree \(m_p\ge2\).
It satisfies the frozen equation
\begin{equation}\label{eq:pos-frozen}
 \sum_{i,j=1}^n B_{ij}(0)\partial_{ij}P_p=0,
\end{equation}
where \(B=A\) in an interior translation chart, and \(B\) is
the transformed matrix in a boundary chart.
\end{lemma}

\begin{proof}
If $p$ is an interior point, if every derivative vanished, then locally \(|\Phi(x)|\le C_N|x-p|^N\) for every \(N\). By strong unique continuation, \(\Phi\)  vanishes on the connected domain, which is a contradiction.
If $p$ is a boundary point, we \textbf{claim} that the vanishing order of $\Phi$ at $p$ is also finite. In details, we assume $\gamma(y^{\prime})$ parameterize the boundary near $p$ with inward unit normal $\nu(y^{\prime})$. For $s>0$, let \[
 \Theta(y',s)=\gamma(y')+sA(\gamma(y'))\nu(y')
\]
 with $\Theta(0)=p$. Its first \(n-1\) columns at \(s=0\) are tangent to the boundary. The last column has positive inward normal component
\(\nu^TA\nu\ge\kappa\). The inverse function theorem gives required chart after shrinking it. Denote \(J=|\det D\Theta|\). The transformed equation is
\[
 {\rm div}_y(B\nabla_y\widehat\Phi)+\lambda_*J\widehat\Phi=0, \;\text{where}\;\widehat\Phi=\Phi\circ\Theta,\;
 B=J(D\Theta)^{-1}A(\Theta)(D\Theta)^{-t}.
\]
Let \(Q=\operatorname{diag}(1,\ldots,1,-1)\). For \(s<0\), define
\[
 \widetilde B(y',s)=Q B(y',-s)Q,\quad
 \widetilde J(y',s)=J(y',-s),\quad
 \widetilde\Phi(y',s)=-\widehat\Phi(y',-s),
\]
Tangential entries and the \((n,n)\) entry are reflected evenly while
mixed entries are reflected oddly. Their zero boundary trace makes the mixed extensions Lipschitz. All other extensions are Lipschitz as well. One can check that $\widetilde{\Phi}$ satisfies 
\begin{equation}
 {\rm div}(\widetilde B\nabla\widetilde\Phi)
       +\lambda_*\widetilde J\widetilde\Phi=0,
\end{equation}
where $\widetilde{B}$ is Lipschitz and uniformly elliptic. Now, also by strong unique continuation, the vanishing order of $\widetilde{\Phi}$ at $0$ is finite, and therefore the vanishing order of $\Phi$ at $p$ is finite, hence the claim holds. Furthermore, a first nonzero Taylor term exists, and we denote it as $P_{p}$, the criticality gives \(m_p\ge2\). Denote $B$ as the transformed matrix, it is easy to see
\begin{equation*}
     \sum_{i,j=1}^n B_{ij}(0)\partial_{ij}P_p=0
\end{equation*}
In addition, let $T=B(0)^{1/2}$, then $\widetilde{P}_{p}(\xi):=P_{p}(T\xi)$ is harmonic.
\end{proof}

Next, we prove a finite-order extension across a flat boundary.
\begin{lemma}
\label{lem:pos-extension}
For every integer \(M\ge0\), $R>0$ and $r\leq R/(2M+1)$. Let $B_{R}^{+}=B_{R}\cap \mathbb{R}_{+}^{n}$, then for any $f\in C^{M}(\overline{B_{R}^{+}})$, there exists a bounded linear extension 
\begin{equation*}
    E_{M}:C^{M}(\overline{B_{R}^{+}})\rightarrow C^{M}(\overline{B_{r}})
\end{equation*}
such that $(1)$ $E_{M}f(y^{\prime},s)=f(y^{\prime},s)$ for $(y^{\prime},s)\in B_{R}^{+}$ and \\
\indent $(2)$ for any multi-index $\beta\in\mathbb{N}^{n-1}$ and integer $k\geq 0$ with $|\beta|+k\leq M$, we have 
\begin{equation*}
\partial_{y^{'}}^{\beta}\partial_{s}^{k}E_{M}f(y^{\prime},0)=\partial_{y^{'}}^{\beta}\partial_{s}^{k}f(y^{\prime},0)
\end{equation*}
\end{lemma}
\begin{proof}
For $s<0$, we define
\begin{equation}
 E_Mf(y',s)=\sum_{j=1}^{M+1}c_jf(y',-js),
\end{equation}
and set \(E_Mf=f\) for \(s\ge0\). Choose the constants by
$ \sum_{j=1}^{M+1}c_j(-j)^k=1, 0\le k\le M$.
It is easy to see the Vandermonde matrix is invertible. Shrink the full neighborhood so that all reflected arguments lie in the original half-neighborhood. For \(|\beta|+k\le M\), the negative-side derivative is
\begin{equation}
 \sum_{j=1}^{M+1}c_j(-j)^k
       \partial_{y'}^\beta\partial_s^kf(y',-js).
\end{equation}
At \(s=0\), it equals the positive-side derivative by the
coefficient identities. Derivatives through order \(M\) therefore
match, and the finite sum also proves the \(C^M\) norm bound.
\end{proof}

\begin{remark}
The odd reflection is used only for unique continuation.
It need not be smooth to high order.
The extension \(E_M\) is used only for the geometric perturbation
lemma. It need not preserve the differential equation.
Separating these two operations avoids an unjustified smooth
reflection assumption.
\end{remark}
\subsection{Singular set comparison of eigenfunction}

\begin{lemma}
\label{lem:pos-global-stability}
There exist an integer \(M_\Phi\ge1\), \(\eta_\Phi\in(0,1)\), and a positive constant \(C_\Phi\) such that if \(v\in C^{M_\Phi}(\overline{\Omega})\) satisfies $||{v-\Phi}||_{C^{M_\Phi}(\overline{\Omega})}<\eta_{\Phi}$, then its global singular set
\[
 K_v=\{x\in\overline{\Omega}:v(x)=0,\ \nabla v(x)=0\}
\]
can be covered by countably many \(C^1\) embedded \((n-2)\)-dimensional submanifolds, and
\begin{equation}\label{eq:pos-global-stable-measure}
 \mathcal{H}^{n-2}(K_v)\le C_\Phi.
\end{equation}
\end{lemma}

\begin{proof}
\textbf{Step 1: prepare the conditions in Lemma \ref{lem:pos-harmonic-stability} locally.}
For \(p\in K_\Phi\), let \(m=m_p\) be the finite order in
Lemma~\ref{lem:pos-finite-order}, and put
\[
 M_p=2(m-1)^2+1\ge m.
\]
Let \(x=\Theta_p(y)\) be a translation chart in the interior
or the corresponding chart at the boundary. In the boundary case, apply \(E_{M_p}\) to the pulled-back function, while in the interior no extension is needed. Then denote \(y=T\xi\), with \(T=B(0)^{1/2}\) and \(\widetilde\Phi_p\) as in Lemma \ref{lem:pos-finite-order}. It is \(C^{M_p}\) on a full neighborhood, and its first homogeneous polynomial \(\widetilde P_p\) is harmonic. The extension also preserves its boundary jet.

 Apply Lemma~\ref{lem:pos-harmonic-stability} to
\(\widetilde P_p\) with \(\varepsilon_p\) and \(r_p\). By rescaled Taylor convergence, there exists a fixed scale \(s_p>0\) such that all coordinate maps are defined and
\[
 ||{s_p^{-m}\widetilde\Phi_p(s_p\,\cdot)
             -\widetilde P_p}||_{C^{M_p}(B_1)}
 <\varepsilon_p/2.
\]
Performing the same operations on \(v\), we obtain \(\widetilde v_p\). Fixed charts, a bounded extension operator, and a fixed scale give 
\[
 ||s_p^{-m}(\widetilde v_p-\widetilde\Phi_p)
                         (s_p\,\cdot)||_{C^{M_p}(B_1)}
 \leq C||{v-\Phi}||_{C^{M_p}(\overline{\Omega})}.
\]
Choose \(\eta_p>0\) with \(C\eta_p<\varepsilon_p/2\).
Then the rescaled function is within \(\varepsilon_p\) of the
harmonic model whenever \(||{v-\Phi}||_{C^{M_p}}<\eta_p\).

\textbf{Step 2: transfer the uniform local measure bound.}
Define the smooth local diffeomorphism
\[
 F_p(\xi)=\Theta_p(s_pT\xi),\qquad
 U_p=F_p(B_{r_p}).
\]
Shrink the original chart, if needed before choosing \(s_p\),
so that \(F_p\) is defined beyond this closed ball.
Its Lipschitz constant on the ball is finite and is
independent of \(v\).

 On the part corresponding to \(\overline{\Omega}\), the extension
equals the original function and matches its boundary derivatives.
The identity
\[
 \nabla_\xi(v\circ F_p)
       =(DF_p)^{t}(\nabla_xv)\circ F_p
\]
and invertibility of \(DF_p\) show that every point of
\(K_v\cap U_p\) comes from a critical point of the rescaled
extended function. This also holds on \(\partial\Omega\).
By the Lipschitz property and
Lemma~\ref{lem:pos-harmonic-stability},
\begin{equation}\label{eq:pos-local-transferred-bound}
 \mathcal{H}^{n-2}(K_v\cap U_p)
 \leq C(n)(m_p-1)^2
                r_p^{n-2}
 =: C_p.
\end{equation}
The graph covering transfers under the same diffeomorphism,
after restricting the covering manifolds to its coordinate domain.
Thus both \(U_p\) and $C_{p}$ are fixed independently
of \(v\).

\textbf{Step 3: obtain the global measure bound.}
If \(K_\Phi\ne\varnothing\), its compactness gives
\[
 K_\Phi\subset U_{p_1}\cup\cdots\cup U_{p_N}.
\]
We choose
\[
 M_\Phi=\max(1,M_{p_1},\ldots,M_{p_N})\;\text{and}\;
 C_\Phi=\sum_{k=1}^N C_{p_{k}}.
\]
On the complement $
 F=\overline{\Omega}\setminus\bigcup_{k=1}^N U_{p_{k}}
$, there is no singularity of \(\Phi\).
If \(F\ne\varnothing\), then
\[
 b=\min_F(|\Phi|+|\nabla\Phi|)>0.
\]
Choose \(\eta_\Phi\) smaller than every \(\eta_{p_{k}}\), and smaller than \(b/2\) when the complement is nonempty.
The sum-of-suprema norm convention ensures on \(F\) that
\[
 |v|+|\nabla v|
 \ge |\Phi|+|\nabla\Phi|-||{v-\Phi}||_{C^1}>b/2.
\]
Hence \(K_v\) lies in the finite union of the \(U_{p_{k}}\).
Therefore,  by subadditivity and \eqref{eq:pos-local-transferred-bound}, we obtain
\eqref{eq:pos-global-stable-measure}, and the local graph
coverings give the asserted global covering.
If \(K_\Phi=\varnothing\), use the positive minimum of
\(|\Phi|+|\nabla\Phi|\) on \(\overline{\Omega}\) directly.
Then \(M_\Phi=1\) and \(C_\Phi=0\) are admissible.
\end{proof}
Note that the above constant $C_{\Phi}$ will depend on the eigenfunction $\Phi$. In the following, we prove the critical set comparison lemma for normalized eigenspace, which helps us to obtain a uniform measure bound.
\begin{lemma}\label{lemma4.7}
\label{lem:pos-eigensphere}
Fix a Dirichlet eigenvalue \(\lambda\) and let
\[
 \mathcal E_\lambda
 =\{\Psi\,\text{real}:L\Psi=\lambda\Psi,\ \Psi|_{\partial\Omega}=0,
                \ ||\Psi||_{L^{2}(\Omega)}=1\}.
\]
There exist an integer \(M_\lambda\), $\eta_{\lambda}\in(0,1)$ and a positive constant \(C_\lambda\),
depending only on \(n,\lambda,A,\Omega\), such that if 
\[
 \inf_{\Psi\in\mathcal E_\lambda}
       ||{v-\Psi}||_{C^{M_\lambda}(\overline{\Omega})}<\eta_\lambda,
\]
then \(K_v\) has the preceding graph covering and \(\mathcal{H}^{n-2}(K_v)\le C_\lambda\).
\end{lemma}

\begin{proof}
The eigenspace is finite-dimensional and consists of smooth
functions. Its unit sphere \(\mathcal E_\lambda\) is compact,
and every finite \(C^M\) norm is continuous on it.
For each \(\Psi\in\mathcal E_\lambda\), take
\(M_\Psi,\eta_\Psi,C_\Psi\) from
Lemma~\ref{lem:pos-global-stability}.
Since  the relative neighborhoods
\[
 \left\{\Psi'\in\mathcal E_\lambda:
  ||{\Psi'-\Psi}||_{C^{M_\Psi}}<\eta_\Psi/2\right\}
\]
cover the sphere. We can choose a finite subcover with centers
\(\Psi_1,\ldots,\Psi_J\), and denote
\[
 M_\lambda=\max_jM_{\Psi_j},\quad
 \eta_\lambda=\frac12\min_j\eta_{\Psi_j},\quad
 C_\lambda=\max_j C_{\Psi_j}.
\]
If the stated distance is less than \(\eta_\lambda\), choose
\(\Psi'\) realizing a distance less than that tolerance $\eta_{\lambda}$,
and then a covering center \(\Psi_j\). Therefore
\[
 ||{v-\Psi_j}||_{C^{M_{\Psi_j}}}
 \le||{v-\Psi'}||_{C^{M_\lambda}}
       +||{\Psi'-\Psi_j}||_{C^{M_{\Psi_j}}}
 <\eta_{\Psi_j}.
\]
Now by the choice of constants, Lemma \ref{lem:pos-global-stability} gives the desired covering and measure bound.
\end{proof}
Now, we are ready to prove Theorem \ref{thm:pos-main}.
\begin{proof}[Proof of Theorem~\ref{thm:pos-main}]
By Lemma~\ref{lem:pos-profile}, we can write
\[
 e^{\lambda_*t}u(t)=\Phi+R(t),\, \text{and denote}\,
 a=||\Phi||_{L^2}(\Omega)>0.
\]
Let 
\[
 v(t)=e^{\lambda_*t}u(t)/a,\qquad
 \Psi=\Phi/a\in\mathcal E_{\lambda_*}.
\]
Lemma~\ref{lem:pos-eigensphere} gives
\(M_{\lambda_*},\eta_{\lambda_*},C_{\lambda_*}\). The high-order spectral estimate yields
\[
 ||{v(t)-\Psi}||_{C^{M_{\lambda_*}}(\overline{\Omega})}
 \le a^{-1}C_{M_{\lambda_*}}e^{-\delta t},\qquad t\ge1.
\]
Choose \(T\ge1\) such that the right-hand side is less than
\(\eta_{\lambda_*}\) for all $t\geq T$. By Lemma \ref{lemma4.7}, we have the graph covering and the uniform bound \(\mathcal{H}^{n-2}(K_{v(t)})\le C_{\lambda_*}\). Since multiplication by the nonzero scalar
\(a^{-1}e^{\lambda_*t}\) does not change critical zeros, then \(K_{v(t)}=\Sigma_t\). Therefore, we have finished the proof of Theorem \ref{thm:pos-main}.
\end{proof}

\part*{\texorpdfstring{\centering\Large Part II. A Counterexample at Low Regularity}{Part II. A Counterexample at Low Regularity}}

Throughout Part II, let \(B_1\) denote the unit ball in
\(\mathbb R^3\). We write \(X=(x,y,z)\),
\(X'=(x,y)^T\), \(r=|X'|\), and \(R=|X|\).
Unless otherwise stated, the objects constructed in each proposition
below are fixed for subsequent use, and their notation is retained
throughout Part II. 

We consider an auxiliary self-adjoint
Dirichlet operator of the form
\begin{equation*}
K_4=-\rho^{-1}\operatorname{div}(\rho A_3\nabla)+V,
\qquad \rho>0.
\end{equation*}
Suppose that \(K_4\) admits two Dirichlet eigenfunctions \(\phi\) and
\(\psi\), with eigenvalues \(\lambda<\mu\):
\begin{equation*}
K_4\phi=\lambda\phi,\qquad K_4\psi=\mu\psi,
\qquad q=\mu-\lambda.
\end{equation*}
The two eigenfunctions yield a solution to the auxiliary Dirichlet
evolution equation \(\partial_t U+K_4U=0\), given by
\begin{equation*}
U(X,t)=e^{-\lambda t}\phi(X)-e^{-\mu t}\psi(X).
\end{equation*}
We seek a fixed operator and eigenfunction pair such that, on infinitely
many disjoint cylinders
\begin{equation*}
\mathcal C_j=\{r<\varepsilon_j,\ |z-z_j|<\varepsilon_j\},
\end{equation*}
we have
\begin{equation*}
\phi(X)=a_j\cosh\bigl(\sqrt q\,(z-z_j)\bigr),
\qquad \psi(X)=1,
\end{equation*}
where \(0<a_j<1\) and \(a_j\to0\). Then, on \(\mathcal C_j\),
\begin{equation*}
U(X,t)=a_j e^{-\lambda t}\cosh\bigl(\sqrt q\,(z-z_j)\bigr)
           -e^{-\mu t}.
\end{equation*}
In particular, at the times \(t_j=q^{-1}\log(1/a_j)\to\infty\),
we obtain
\begin{equation*}
U(\cdot,t_j)=0,\qquad \nabla U(\cdot,t_j)=0
\quad\text{on } D_j=\{z=z_j,\ r<\varepsilon_j\}.
\end{equation*}
The auxiliary solution \(U\) already has the essential features of
the desired counterexample. Finally, we modify \(U\) to obtain
the counterexample \(u\).

Part II is organized as follows. In Section 5, modifying the classical eigenfunctions, we obtain an elliptic operator with a potential and a density, and its eigenfunction pair is $(\phi_{2},\psi_{2})=(P(x,y),1)$. As mentioned before, $\psi=1$ need not change, while we hope $\phi$ has the form ${\rm cosh(z)}$ in local cylinder. Hence Section 6 prepares the eigenfunctions, which locally involve only $P$ and ${\rm cosh(z)}$. In Section 7, we remove the transverse polynomial $P(x,y)$ by cut-off and correction, and obtain eigenfunctions involving only $a\,{\rm cosh(z)}$. We then show that the corresponding matrix is uniformly elliptic and $C^{1,\gamma}$ for any $\gamma\in(0,1/2)$, but lacks higher regularity.  Section 8 removes the potential and density in elliptic operator and complete the construction. Ii is worth mentioning that, since the construction is relatively complicated, \textbf{we summarize what this step entails before each key step and highlight it in italics.}

\section{A Dirichlet eigenfunction pair with an explicit local matrix}\label{construction-of-the-counterexample-at-low-regularity}

% --------------------------------------------------------------------
% Preparing a global Dirichlet eigenfunction pair
\subsection{Preparing a global Dirichlet eigenfunction pair}\label{preparing-a-global-dirichlet-eigenfunction-pair}

We first record the Dirichlet eigenfunctions used in the construction.

\begin{proposition}\label{prop:initial-dirichlet-pair}
There exist two Dirichlet eigenfunctions \(\phi_0,\psi_0\) of
\(-\Delta\) on \(B_1\), with respective eigenvalues
\(\lambda<\mu=16\pi^2\), such that
\begin{equation*}
\phi_0(X)=a(R)P(x,y)
\qquad\text{and}\qquad
\psi_0(X)=\frac{\sin(4\pi R)}{4\pi R},
\end{equation*}
where
\begin{equation*}
P(x,y)=x^4-6x^2y^2+y^4,
\end{equation*}
and \(a(0)=1\). The value \(\psi_0(0)=1\) is understood
by continuity.
\end{proposition}

\begin{proof}
Since \(P\) is a homogeneous harmonic polynomial of degree \(4\),
\begin{equation*}
\Delta(a(R)P)
=\left(a''(R)+\frac{10}{R}a'(R)\right)P.
\end{equation*}
Let \(\lambda\) be the first Dirichlet eigenvalue of the unit ball
in \(\mathbb R^{11}\). Its positive first eigenfunction is radial,
by rotational invariance and simplicity of the first eigenvalue.
Denote its radial profile by \(a(R)\), normalized so that \(a(0)=1\).
By elliptic regularity, \(a\) is smooth at the origin, and the
radial eigenvalue equation is
\begin{equation*}
 a''(R)+\frac{10}{R}a'(R)+\lambda a(R)=0,
 \qquad a'(0)=0,\quad a(1)=0.
\end{equation*}
Thus \(\phi_0(X)=a(R)P(x,y)\) is a Dirichlet eigenfunction
of \(-\Delta\) on \(B_1\subset\mathbb R^3\).

By the variational characterization of \(\lambda\), taking
\(f(R)=1-R^2\) as a test function gives
\begin{equation*}
\lambda\le\frac{\int_0^1 R^{10}|f'(R)|^2\,dR}
{\int_0^1 R^{10}|f(R)|^2\,dR}=\frac{165}{2}<16\pi^2.
\end{equation*}
The stated function \(\psi_0\) satisfies
\(-\Delta\psi_0=16\pi^2\psi_0\) and vanishes on
\(\partial B_1\). Thus \(\mu=16\pi^2\) and
\(\lambda<\mu\), as required.
\end{proof}

\begin{proposition}\label{prop:global-dirichlet-normal-form}
There exist smooth functions \(\rho>0,V\) and a smooth symmetric
uniformly elliptic matrix \(A_0\) on \(B_1\) such that the Dirichlet
operator
\begin{equation*}
K_1=-\rho^{-1}\operatorname{div}(\rho A_0\nabla)+V
\end{equation*}
admits eigenfunctions \(\phi_2,\psi_2\) with eigenvalues
\(\lambda,\mu\), respectively. Near the origin,
\begin{equation*}
\rho=1,\quad V=\mu,\quad \phi_2=P,\quad\psi_2=1,\quad
A_0=I+O(|X|^2).
\end{equation*}
Moreover, \(A_0\) is invariant under rotations and reflections in
the \((x,y)\)-plane.
\end{proposition}

\begin{proof}
\textbf{\textit{Conjugation.}}
Choose a smooth positive radial function \(h_0\) on \(\overline{B_1}\)
equal to \(\psi_0\) near the origin, and set
\begin{equation*}
\phi_1=\frac{\phi_0}{h_0},\qquad \psi_1=\frac{\psi_0}{h_0}.
\end{equation*}
Then \(K_0\phi_1=\lambda\phi_1\) and \(K_0\psi_1=\mu\psi_1\), where
\begin{equation*}
K_0=h_0^{-1}(-\Delta)h_0
=-h_0^{-2}\operatorname{div}(h_0^2\nabla)+V_0,
\qquad V_0=-\frac{\Delta h_0}{h_0}.
\end{equation*}
Near the origin, \(\phi_1=(a(R)/h_0(X))P\),
\(\psi_1=1\), and \(V_0=\mu\).

\textbf{\textit{Local normalization.}}
First, set
\begin{equation*}
(\xi_1,\xi_2)
=\left(\frac{a(R)}{h_0(X)}\right)^{1/4}(x,y),
\end{equation*}
leaving \(z\) unchanged. Since \(P\) is homogeneous of degree \(4\),
\begin{equation*}
P(\xi_1,\xi_2)=\frac{a(R)}{h_0(X)}P(x,y)=\phi_1(X).
\end{equation*}
The density after this change is
\begin{equation*}
\rho_1(\xi_1,\xi_2,z)
=\frac{h_0(X)^2}{\displaystyle
\det\frac{\partial(\xi_1,\xi_2,z)}{\partial(x,y,z)}}.
\end{equation*}
Next, set
\begin{equation*}
\xi_3=\int_0^z\rho_1(\xi_1,\xi_2,s)\,ds.
\end{equation*}
The second change has Jacobian \(\rho_1\), so the combined change
has Jacobian \(h_0^2\). Thus the new density is \(1\), and the
eigenfunctions become \(P\) and \(1\) near the origin.

\textbf{\textit{Global extension and transformed coefficients.}}
Since \(a\) and \(h_0\) are smooth radial functions with
\(a(0)=h_0(0)=1\), we have \(a'(0)=h_0'(0)=0\), and hence
\(a(R)/h_0(X)=1+O(R^2)\). The first change has Jacobian
\(1+O(R^2)\), so \(\rho_1=1+O(R^2)\).
Denoting the combined coordinate change by \(G\), we obtain
\begin{equation*}
G(X)=X+O(|X|^3),\qquad DG(X)=I+O(|X|^2).
\end{equation*}
Using a cutoff, we extend \(G\) to a smooth diffeomorphism of
\(B_1\) that equals the identity near the boundary.

Define
\begin{equation*}
\phi_2=\phi_1\circ G^{-1},\qquad
\psi_2=\psi_1\circ G^{-1}.
\end{equation*}
The transformed coefficients are
\begin{equation*}
A_0=(DG\,DG^T)\circ G^{-1},\qquad
\rho=\left(\frac{h_0^2}{\det DG}\right)\circ G^{-1},
\qquad V=V_0\circ G^{-1}.
\end{equation*}
The required local properties follow from these formulas.
The radial choice of \(h_0\) and the radial cutoff preserve rotations
and reflections in the \((x,y)\)-plane, so \(A_0\) has the stated
symmetries.
\end{proof}

The second eigenfunction $\psi_2$ already has the local form needed for
the construction. We keep it unchanged while modifying the first
eigenfunction and the coefficient matrix, until the final transformations
in Section~\ref{removing-the-potential-and-density}.

% --------------------------------------------------------------------
% An explicit local principal matrix and compatible patching
\subsection{An explicit local matrix preserving the eigenfunctions}\label{an-explicit-local-principal-matrix-and-compatible-patching}

We modify \(A_0\) near the origin to obtain an explicit coefficient
matrix suitable for the subsequent construction, while preserving
the two Dirichlet eigenfunctions and their eigenvalues.

Recall that \(X'=(x,y)^T\), and define
\begin{equation*}
A_*=
\begin{pmatrix}
I_2+\frac q{24}X'(X')^T&0\\
0&1
\end{pmatrix}.
\end{equation*}

\begin{proposition}\label{prop:explicit-local-matrix}
For every sufficiently small \(\delta>0\), there exists a smooth,
symmetric, uniformly elliptic matrix \(A_1\) on \(B_1\) such that
\begin{equation*}
A_1=
\begin{cases}
A_*, & r<\delta,\ |z|<\delta,\\
A_0, & r\ge 2\delta\ \text{or}\ |z|\ge 2\delta.
\end{cases}
\end{equation*}
Moreover, \(\phi_2,\psi_2\) remain Dirichlet eigenfunctions of
\[
K_2=-\rho^{-1}\operatorname{div}(\rho A_1\nabla)+V
\]
with eigenvalues \(\lambda,\mu\), respectively.
\end{proposition}

\begin{proof}
\textbf{\textit{Compatibility of the local matrices.}}
Since $\Delta_{x,y}P=0$ and $X'\cdot\nabla P=4P$, we have
\begin{equation*}
\operatorname{div}(A_*\nabla P)=qP.
\end{equation*}
The same identity holds for $A_0$ near the origin by
Proposition~\ref{prop:global-dirichlet-normal-form}. Thus
\begin{equation*}
D=A_*-A_0,\qquad \operatorname{div}(D\nabla P)=0
\end{equation*}
near the origin. We cut off $D$ while preserving this identity.

\textbf{\textit{Cutoff in the $z$-direction.}}
First, consider the $z$-direction. The matrices $A_0$ and $A_*$
are invariant under rotations and reflections in the $(x,y)$-plane.
For the mixed block $w=(D_{xz},D_{yz})^T$, these symmetries imply
that $w$ is parallel to $X'$, with a proportionality factor depending
only on $r$ and $z$. Smoothness at $r=0$ then yields
\begin{equation*}
w=d(r^2,z)X'
\end{equation*}
for some smooth function $d$. Choose a smooth cutoff $\eta(z)$ equal to $1$
for $|z|<\delta$ and to $0$ for $|z|\ge 2\delta$. Then
\begin{equation*}
\operatorname{div}(\eta D\nabla P)
=\eta' w\cdot\nabla P=4\eta'dP.
\end{equation*}
Define
\begin{equation*}
c(r,z)=-\eta'(z)\int_0^1 t^5d(t^2r^2,z)\,dt,
\qquad
C_z=\begin{pmatrix}cX'(X')^T&0\\0&0\end{pmatrix}.
\end{equation*}
The integral expression is smooth at $r=0$ and gives
\begin{equation*}
6c+r\partial_rc=-\eta'd,\qquad
\operatorname{div}(C_z\nabla P)
=4(6c+r\partial_rc)P=-4\eta'dP.
\end{equation*}
Consequently,
\begin{equation*}
D^{(z)}=\eta D+C_z,
\qquad \operatorname{div}(D^{(z)}\nabla P)=0.
\end{equation*}
Moreover, $D^{(z)}=D$ for $|z|<\delta$ and vanishes for
$|z|\ge2\delta$.

\textbf{\textit{Transverse cutoff and correction.}}
Next, choose a smooth radial cutoff $\chi(r)$ equal to $1$ for
$r<\delta$ and to $0$ for $r\ge2\delta$. Let $\theta$ be the polar
angle. The upper-left $2\times2$ submatrix $B$ of $D^{(z)}$,
corresponding to the $x,y$ coordinates, has the form
$B=\alpha(r,z)I_2+\beta(r,z)X'(X')^T$ by rotational and reflection
symmetry. Since $P=r^4\cos(4\theta)$ is independent of $z$,
\begin{align*}
E:=\operatorname{div}(\chi D^{(z)}\nabla P)
&=\chi'(r)\frac{X'}r\cdot B\nabla_{X'}P\\
&=4r^3\chi'(r)(\alpha+\beta r^2)\cos(4\theta)
=:f(r,z)\cos(4\theta).
\end{align*}
The function $f$ is supported in
$\delta<r<2\delta$, $|z|<2\delta$. Set
\begin{equation*}
J=-\frac r4 f(r,z)\sin(4\theta)e_\theta,
\qquad e_\theta=\frac{(-y,x,0)}r.
\end{equation*}
Then $\operatorname{div}J=-E$. On this annulus,
$v=\nabla P$ satisfies $|v|=4r^3>0$. The symmetric matrix
\begin{equation}\label{eq:symmetric-matrix-correction}
C_r=\frac{Jv^T+vJ^T}{|v|^2}
    -\frac{J\cdot v}{|v|^4}vv^T
\end{equation}
satisfies $C_r\nabla P=J$. Hence, with
$D^{(r)}=\chi D^{(z)}+C_r$, ${\rm div}(D^{(r)}\nabla P)=0$, we may set
$A_1=A_0+D^{(r)}$ and 
extend $D^{(r)}$ by zero outside $\{r<2\delta,\ |z|<2\delta\}$.
The cutoff formulas give $A_1=A_*$ for $r<\delta$, $|z|<\delta$,
and $A_1=A_0$ for $r\ge2\delta$ or $|z|\ge2\delta$.

\textbf{\textit{Ellipticity and the eigenfunction equations.}}
Since $D=O(R^2)$ and $d(0,0)=0$, the cutoff corrections
$C_z,C_r$ are $O(\delta^2)$ on their supports. Thus $A_1$ is
uniformly elliptic for sufficiently small $\delta$.
The perturbation $A_1-A_0$ is supported where $\rho=1$,
$V=\mu$, $\phi_2=P$, and $\psi_2=1$.
Since $\operatorname{div}((A_1-A_0)\nabla P)=0$ and
$\nabla\psi_2=0$ there, both eigenfunction equations remain valid.
\end{proof}
% --------------------------------------------------------------------
% Smooth perturbations with prescribed critical values
\section{Construction of eigenfunctions with prescribed critical values}\label{sec:critical-value-perturbations}
In this section, we modify the eigenfunctions and the coefficient matrix
constructed in Section~5, while preserving the Dirichlet eigenvalue equations.
We obtain infinitely many critical points $X_j=(0,0,z_j)$, where $z_j\to0$,
with positive critical values $a_j\to0$. Near each $X_j$, the modified
eigenfunction has a prescribed local form that will be used to construct
critical disks in the next section.
\subsection{Construction of perturbations in the $z$-direction}\label{signed-controls-and-critical-values-approaching-zero}
The function $f_0=P+a\cosh(\sqrt q\,z)$ satisfies
$\operatorname{div}(A_*\nabla f_0)=qf_0$, with
$f_0(0)=a$ and $\nabla f_0(0)=0$.
However, replacing the added term by a smooth compactly supported
function $h(z)$ while keeping $A_*$ fixed would require $h''=qh$,
which has no nonzero compactly supported solution. We therefore add
quadratic terms and modify the coefficient matrix.

\begin{proposition}\label{prop:s56-controls}
\begin{samepage}
For every sufficiently small $\ell>0$ and every $a>0$, there exist
$g,k,h\in C_c^\infty((-\ell,\ell))$ such that, setting
\[
f_1=P+g(z)xy+k(z)(x^2-y^2)+h(z),
\]
and
\[
A_{**}=A_*+
\begin{pmatrix}
c(z)&b(z)&0\\
b(z)&-c(z)&0\\
0&0&0
\end{pmatrix},
\]
where
\begin{equation}\label{3.0}
b(z)=\frac{3g''(z)-2qg(z)}{144},\qquad
c(z)=\frac{2qk(z)-3k''(z)}{72},
\end{equation}
we have
\begin{equation}\label{eq:s56-local-eigenfunction}
\operatorname{div}(A_{**}\nabla f_1)=qf_1.
\end{equation}
\end{samepage}
Moreover, for $|z|<2\ell/5$,
\begin{equation}\label{3.4}
f_1=P+a\cosh(\sqrt q\,z),\qquad A_{**}=A_*.
\end{equation}
For every integer $m\ge0$,
\begin{equation}\label{eq:s56-control-estimates}
\|g\|_{C^m(\mathbb R)}+\|k\|_{C^m(\mathbb R)}
\le C_m\sqrt a\,\ell^{-m},\qquad
\|h\|_{C^m(\mathbb R)}\le C_ma\,\ell^{-m},
\end{equation}
where the constants are independent of $a,\ell$.
\end{proposition}

\begin{proof}
\textbf{\textit{Reduction to a scalar equation.}}
A direct calculation shows that \eqref{eq:s56-local-eigenfunction} is equivalent to
\begin{equation}\label{3.1}
h''-qh=-\frac{g(3g''-2qg)}{72}
         +\frac{k(3k''-2qk)}{18}.
\end{equation}
Set
\[
H=h+\frac{g^2}{48}-\frac{k^2}{12},\qquad
\mathcal T(w)=\frac{(w')^2}{24}+\frac{qw^2}{144}.
\]
Then \eqref{3.1} becomes
\begin{equation}\label{3.2}
H''-qH=\mathcal T(g)-4\mathcal T(k).
\end{equation}

\textbf{\textit{Moment conditions and compact support.}}
Fix a sufficiently small $\ell>0$ and work on the normalized
interval $(-1,1)$, with coordinate $\tau$. Choose nonzero even functions
$W_+,W_-\in C_c^\infty((-1,1))$ with supports contained, respectively, in
\[
\{7/10<|\tau|<9/10\},\qquad
\{2/5<|\tau|<3/5\}.
\]
Set
\[
F_\pm=\frac{(W_\pm')^2}{24}
       +\frac{q\ell^2W_\pm^2}{144}.
\]
Define
\[
\begin{aligned}
A_\pm&=\int_0^1\cosh(\sqrt q\,\ell t)F_\pm(t)\,dt,\\
B_\pm&=\int_0^1\frac{\sinh(\sqrt q\,\ell t)}{\sqrt q\,\ell}
F_\pm(t)\,dt.
\end{aligned}
\]
\begin{samepage}
Since $t\mapsto\tanh(\sqrt q\,\ell t)/(\sqrt q\,\ell)$ is
strictly increasing and the supports of $F_+$ and $F_-$ are
separated, we have
\[
\frac{B_+}{A_+}>\frac{B_-}{A_-}.
\]
The difference stays bounded away from zero as $\ell\to0$,
while $A_\pm$ have positive limits. Hence there are positive
coefficients $d_+,d_-$, uniformly bounded for small $\ell$, such that
\end{samepage}
\begin{equation}\label{eq:s56-moment-conditions}
d_+A_+-d_-A_-=0,\qquad d_+B_+-d_-B_-=1.
\end{equation}
Set
\[
S=d_+F_+-d_-F_-,\qquad
\widetilde g=\sqrt{d_+}\,W_+,\qquad
\widetilde k=\frac{\sqrt{d_-}}{2}W_-,
\]
and solve
\[
\widetilde H''-q\ell^2\widetilde H=S,
\qquad \widetilde H(0)=1,\quad \widetilde H'(0)=0.
\]
Both $S$ and $\widetilde H$ are even. Solving the ODE gives, for
$\tau\ge9/10$,
\[
\begin{aligned}
\widetilde H(\tau)={}&
\left(1-\int_0^1
\frac{\sinh(\sqrt q\,\ell t)}{\sqrt q\,\ell}S(t)\,dt\right)
\cosh(\sqrt q\,\ell\tau)\\
&+\left(\int_0^1\cosh(\sqrt q\,\ell t)S(t)\,dt\right)
\frac{\sinh(\sqrt q\,\ell\tau)}{\sqrt q\,\ell}.
\end{aligned}
\]
Both coefficients vanish by \eqref{eq:s56-moment-conditions}.
Since $\widetilde H$ is even, it follows that $\widetilde H=0$
for $|\tau|\ge9/10$. Since $S=0$ on $|\tau|<2/5$,
\[
\widetilde H(\tau)=\cosh(\sqrt q\,\ell\tau)
\qquad (|\tau|<2/5).
\]
Define
\[
\widetilde h=\widetilde H
-\frac{\widetilde g^2}{48}
+\frac{\widetilde k^2}{12}.
\]
The functions $\widetilde g,\widetilde k,\widetilde h$ are smooth
and compactly supported in $(-1,1)$.
The uniform bounds on $d_+,d_-$ and the fixed-interval ODE
give uniform $C^m$ bounds for $\widetilde g,\widetilde k,\widetilde h$, for every $m\ge0$.

\textbf{\textit{Rescaling and estimates.}}
Finally, set
\[
g(z)=\sqrt a\,\widetilde g(z/\ell),\qquad
k(z)=\sqrt a\,\widetilde k(z/\ell),\qquad
h(z)=a\,\widetilde h(z/\ell).
\]
Rescaling the normalized equation proves \eqref{3.2}, hence
\eqref{3.1}. The central form \eqref{3.4} and the estimates
\eqref{eq:s56-control-estimates} follow from the construction and
scaling.
\end{proof}

We next apply Proposition~\ref{prop:s56-controls} on disjoint intervals
along the $z$-axis that accumulate at $z=0$, choosing rapidly decaying
amplitudes to retain smoothness.

\begin{corollary}\label{cor:s56-infinite-perturbations}
Let $z_j=2^{-j}$, $\ell_j=2^{-j-4}$, $a_j=e^{-j^4}$,
and $I_j=(z_j-\ell_j,z_j+\ell_j)$.
For every sufficiently large integer $j_0$, there exist a smooth function
$f_2$ and a smooth symmetric, uniformly elliptic matrix $A_{**}$ on a fixed
cylinder about the origin such that
\[
\operatorname{div}(A_{**}\nabla f_2)=qf_2.
\]
For every $j\ge j_0$,
\[
f_2=P+a_j\cosh\bigl(\sqrt q\,(z-z_j)\bigr),\qquad A_{**}=A_*
\quad\text{if }|z-z_j|<2\ell_j/5.
\]
Moreover,
\[
f_2=P,\qquad A_{**}=A_*
\quad\text{if }z\notin\bigcup_{j\ge j_0}I_j.
\]
For every integer $m\ge0$, on any fixed bounded cylinder,
\begin{equation}\label{eq:s56-tail-smallness}
\|f_2-P\|_{C^m}+\|A_{**}-A_*\|_{C^m}
\longrightarrow0\qquad\text{as }j_0\to\infty.
\end{equation}
\end{corollary}

\begin{proof}
Apply Proposition~\ref{prop:s56-controls}, translated by $z_j$, to obtain
$g_j,k_j,h_j\in C_c^\infty(I_j)$ and set
\[
g=\sum_{j\ge j_0}g_j,\qquad
k=\sum_{j\ge j_0}k_j,\qquad
h=\sum_{j\ge j_0}h_j.
\]
The intervals $I_j$ are disjoint, and for sufficiently large $j_0$
all lie in the fixed neighborhood where $A_1=A_*$, $\rho=1$,
$V=\mu$, and $\phi_2=P$.
For every $m,N\ge0$, Proposition~\ref{prop:s56-controls} gives
\[
z_j^{-N}\bigl(\|g_j^{(m)}\|_\infty+\|k_j^{(m)}\|_\infty\bigr)
\le C_{m,N}\exp\!\left(-\frac{j^4}{2}+(m+N)j\log2\right)
\longrightarrow0.
\]
The corresponding bound for $h_j$ has $-j^4$ in place of $-j^4/2$.
These estimates imply that $g,k,h$ extend smoothly across $z=0$,
with all derivatives vanishing there. They also tend to zero in
all $C^m$ norms as $j_0\to\infty$.

Define
\[
f_2=P+g(z)xy+k(z)(x^2-y^2)+h(z)
\]
and let $A_{**}$ have the form in Proposition~\ref{prop:s56-controls},
with $b,c$ given by \eqref{3.0} using the present $g,k$.
The disjoint supports and \eqref{eq:s56-local-eigenfunction} give
$\operatorname{div}(A_{**}\nabla f_2)=qf_2$, first away from $z=0$
and then at $z=0$ by smoothness. Since the intervals $I_j$ are disjoint, the stated formulas in the
central regions and outside their union follow directly from
Proposition~\ref{prop:s56-controls}.
Since $b,c$ involve at most two derivatives of $g,k$,
the preceding estimates give \eqref{eq:s56-tail-smallness}.
Finally, $A_{**}-A_*$ is small in $C^0$ as $j_0\to\infty$;
hence $A_{**}$ is uniformly elliptic on the fixed cylinder for
sufficiently large $j_0$.
\end{proof}

\subsection{Cutoffs and correction of the coefficients}\label{correcting-a-zero-integral-error-on-a-fixed-annulus}

The perturbations in Corollary~\ref{cor:s56-infinite-perturbations}
are localized only in the $z$-direction. We now introduce cutoffs in the
$(x,y)$-variables and correct the resulting error by modifying the
coefficient matrix, obtaining a global Dirichlet eigenfunction.

\begin{proposition}\label{prop:s56-global-cutoff}
Let $z_j,\ell_j,a_j$ be as in
Corollary~\ref{cor:s56-infinite-perturbations}.
For every sufficiently large integer $j_0$, there exist a smooth,
symmetric, uniformly elliptic matrix $A_2$ on $B_1$ and a smooth
function $\phi_3$ such that $\phi_3$ and $\psi_2$ are Dirichlet
eigenfunctions of
\[
K_3=-\rho^{-1}\operatorname{div}(\rho A_2\nabla)+V
\]
with eigenvalues $\lambda,\mu$, respectively.
Outside a fixed small cylinder about the origin,
$A_2=A_1$ and $\phi_3=\phi_2$.
Moreover, for some fixed $r_0>0$ and every $j\ge j_0$,
\[
A_2=A_*,\qquad
\phi_3=P+a_j\cosh\bigl(\sqrt q\,(z-z_j)\bigr)
\]
whenever $r<r_0$ and $|z-z_j|<2\ell_j/5$.
\end{proposition}

\begin{proof}
\textbf{\textit{Transverse cutoffs.}}
Let $g,k,h$ be as in the proof of
Corollary~\ref{cor:s56-infinite-perturbations}, and define
$b,c$ by \eqref{3.0}. Fix $\delta>0$ so that
$A_1=A_*$, $\rho=1$, $V=\mu$, $\phi_2=P$, and $\psi_2=1$
in $\{r<\delta,\ |z|<\delta\}$.
For sufficiently large $j_0$, all $I_j$ lie in $(-\delta,\delta)$.
Choose smooth radial cutoffs $\chi,\zeta$, equal to $1$ near $r=0$
and supported in $r<\delta$, with $\zeta$ signed and
\[
\int_{\mathbb R^2}\zeta(r)\,dx\,dy=0.
\]
In this cylinder, set
\[
f_3=P+\chi(r)\bigl(g(z)xy+k(z)(x^2-y^2)\bigr)+\zeta(r)h(z).
\]
Define the intermediate matrix globally by
\[
\widetilde A_1=A_1+\chi(r)
\begin{pmatrix}
c(z)&b(z)&0\\
b(z)&-c(z)&0\\
0&0&0
\end{pmatrix}.
\]
Here $(\widetilde A_1)_{zz}=1$ and
$(\widetilde A_1)_{zx}=(\widetilde A_1)_{zy}=0$
in $\{r<\delta,\ |z|<\delta\}$. Where $\chi=\zeta=1$,
we have $f_3=f_2$ and $\widetilde A_1=A_{**}$, as in
Corollary~\ref{cor:s56-infinite-perturbations}.
Where both cutoffs vanish, $f_3=P$ and $\widetilde A_1=A_1$.

\textbf{\textit{Zero integral of the error.}}
The cutoff error
\[
E=\operatorname{div}(\widetilde A_1\nabla f_3)-qf_3
\]
is supported in $\{r_1\le r\le r_2\}$ for some fixed
$0<r_1<r_2<\delta$.
Recall that $\operatorname{div}(A_*\nabla P)=qP$.
Subtracting this equation gives
\[
E=\operatorname{div}\bigl(\widetilde A_1\nabla f_3-A_*\nabla P\bigr)-q(f_3-P).
\]
The $z$-component of the flux difference is $\partial_z(f_3-P)$,
and its transverse components are compactly supported in the
$(x,y)$-plane. Thus, integrating over $\mathbb R^2$ yields
\[
\int_{\mathbb R^2}E\,dx\,dy
=(\partial_z^2-q)\int_{\mathbb R^2}(f_3-P)\,dx\,dy.
\]
Because $\chi$ is radial, the angular integrals of $xy$ and
$x^2-y^2$ vanish. By the choice of $\zeta$,
\[
\int_{\mathbb R^2}(f_3-P)\,dx\,dy
=h(z)\int_{\mathbb R^2}\zeta(r)\,dx\,dy=0.
\]
Consequently, $\int_{\mathbb R^2}E(x,y,z)\,dx\,dy=0$ for every $z$.

\textbf{\textit{Flux and matrix correction.}}
Define the angular average of $E$ by
\[
E_0(r,z)=\frac1{2\pi}\int_0^{2\pi}
E(r\cos\theta,r\sin\theta,z)\,d\theta.
\]
In polar coordinates, set
\[
J_r=-\frac1r\int_0^r sE_0(s,z)\,ds,\qquad J_z=0,
\]
and
\[
J_\theta(r,\theta,z)
=-r\int_0^\theta\bigl(E(r\cos\varphi,r\sin\varphi,z)-E_0(r,z)\bigr)\,d\varphi.
\]
Since $E-E_0$ has zero angular mean, $J_\theta$ is $2\pi$-periodic.
The polar-coordinate divergence formula gives
\[
\operatorname{div}J
=\frac1r\partial_r(rJ_r)+\frac1r\partial_\theta J_\theta
=-E_0-(E-E_0)=-E.
\]
Since $E$ is supported in a fixed annulus and has zero integral
on each transverse plane, $J$ is smooth and supported in that annulus.

Since the cutoffs are fixed, Corollary~\ref{cor:s56-infinite-perturbations}
implies that $\nabla_{x,y}f_3$ converges uniformly to
$\nabla_{x,y}P$ as $j_0\to\infty$.
As $|\nabla_{x,y}P|=4r^3$ is bounded below by a positive
constant on the fixed annulus, $\nabla_{x,y}f_3$ is nonzero
there for sufficiently large $j_0$.

Applying the symmetric matrix formula
\eqref{eq:symmetric-matrix-correction} in the $(x,y)$-coordinates,
with $v=\nabla_{x,y}f_3$, we obtain a smooth symmetric correction matrix $C$,
supported in the annulus and having only an upper-left $2\times2$
block, such that $C\nabla f_3=J$. Consequently,
\[
\operatorname{div}\bigl((\widetilde A_1+C)\nabla f_3\bigr)=qf_3.
\]

\textbf{\textit{Global eigenfunctions and ellipticity.}}
Set $A_2=\widetilde A_1+C$ and $\phi_3=\phi_2+(f_3-P)$,
extending the local modifications by zero outside the cylinder.
The corrections tend to zero in $C^0$, so $A_2$ remains uniformly
elliptic. Since $\rho=1$, $V=\mu$, and $\psi_2=1$ on their support,
the eigenvalue equations for $\phi_3$ and $\psi_2$ follow.
The Dirichlet boundary conditions are unchanged, and the prescribed
local forms follow from Corollary~\ref{cor:s56-infinite-perturbations}.
\end{proof}

% --------------------------------------------------------------------
% Turning the critical set of a critical value into a disk
\section{Construction of infinitely many critical disks}\label{turning-the-critical-set-of-a-critical-value-into-a-disk}

We first construct a smooth local modification producing one critical disk.
We then combine infinitely many such modifications and examine the
regularity of the resulting coefficient matrix and eigenfunction.

\subsection{A local construction}

Recall the polynomial $P$ and the matrix $A_*$ from
Section~\ref{construction-of-the-counterexample-at-low-regularity}.
Near the $j$-th center, with $s=z-z_j$ and $a=a_j>0$,
the preceding construction gives
\[
\phi_3=P(x,y)+a\cosh(\sqrt q\,s),\qquad
\operatorname{div}(A_*\nabla\phi_3)=q\phi_3,
\]
and $\psi_2=1$. We now remove the transverse polynomial $P$
on a small cylinder and correct the equation by modifying only
the $(3,3)$ entry of the coefficient matrix.

\begin{proposition}\label{prop:s78-one-patch}
There exists a constant $\delta_*>0$ such that, for every $a>0$ and
$0<\varepsilon\le1$ with $\varepsilon^2/a\le\delta_*$, there are
smooth functions $f_{\varepsilon,a}$ and $d_{\varepsilon,a}$ on
$\mathbb R^2\times\mathbb R$ satisfying
\[
f_{\varepsilon,a}=P+a\cosh(\sqrt q\,s),\qquad d_{\varepsilon,a}=0
\quad\text{outside }\{|X'|<2\varepsilon,\ |s|<2\varepsilon\},
\]
and
\[
f_{\varepsilon,a}=a\cosh(\sqrt q\,s),\qquad
 d_{\varepsilon,a}=0
\quad\text{on }\{|X'|<\varepsilon,\ |s|<\varepsilon\}.
\]
In particular, $f_{\varepsilon,a}=a$ and
$\nabla f_{\varepsilon,a}=0$ on the disk
$\{|X'|<\varepsilon,\ s=0\}$.
The symmetric matrix
\[
A_{\varepsilon,a}=A_*+
\begin{pmatrix}
0&0&0\\
0&0&0\\
0&0&d_{\varepsilon,a}
\end{pmatrix}
\]
satisfies
\[
\operatorname{div}(A_{\varepsilon,a}\nabla f_{\varepsilon,a})
=qf_{\varepsilon,a}.
\]
For every multi-index $\alpha\in\mathbb N_0^3$,
\begin{equation}\label{eq:s78-d-estimate}
\|\partial^\alpha d_{\varepsilon,a}\|_\infty
\le C_\alpha\frac{\varepsilon^{2-|\alpha|}}a,
\end{equation}
where $\partial^\alpha$ acts on $(x,y,s)$ and the constants
$C_\alpha$ are independent of $a$ and $\varepsilon$.
\end{proposition}

\begin{proof}
\textbf{\textit{Construction of the local function.}}
Choose a radial $\chi\in C_c^\infty(\mathbb R^2)$ equal to $1$
on $|X'|\le1$ and supported in $|X'|<2$, and an even
$\eta\in C_c^\infty((-2,2))$, equal to $1$ on $[-1,1]$,
such that
\begin{equation}\label{eq:s78-eta}
\int_0^\infty\eta(t)\,dt=0.
\end{equation}
The latter is obtained by subtracting a suitable outer bump from
an even cutoff equal to $1$ on $[-1,1]$. Set
\[
W_\varepsilon(X')=\chi(X'/\varepsilon)P(X'),\qquad
\eta_\varepsilon(s)=\eta(s/\varepsilon),
\]
\begin{equation}\label{eq:s78-flattened}
f_{\varepsilon,a}(X',s)=P(X')+a\cosh(\sqrt q\,s)
                  -W_\varepsilon(X')\eta_\varepsilon(s).
\end{equation}
\textbf{\textit{Correction of the equation.}}
Let
\begin{equation}\label{eq:s78-K}
\mathcal K_\varepsilon
=\operatorname{div}_{x,y}((A_*)_{xy}\nabla W_\varepsilon)
 -qW_\varepsilon.
\end{equation}
Since $P+a\cosh(\sqrt q\,s)$ solves the unmodified equation,
\begin{equation}\label{eq:s78-error}
E:=\operatorname{div}(A_*\nabla f_{\varepsilon,a})-qf_{\varepsilon,a}
=-\mathcal K_\varepsilon\eta_\varepsilon
 -W_\varepsilon\eta_\varepsilon''.
\end{equation}
Define
\[
I_\varepsilon(s)=\int_0^s\eta_\varepsilon(t)\,dt,\qquad
N=\mathcal K_\varepsilon I_\varepsilon
   +W_\varepsilon\eta_\varepsilon'.
\]
By the definition of $A_{\varepsilon,a}$,
\[
\operatorname{div}(A_{\varepsilon,a}\nabla f_{\varepsilon,a})
=\operatorname{div}(A_*\nabla f_{\varepsilon,a})
 +\partial_s(d_{\varepsilon,a}\partial_s f_{\varepsilon,a}).
\]
Thus, to eliminate $E$, it suffices to require
\[
\partial_s(d_{\varepsilon,a}\partial_s f_{\varepsilon,a})=-E.
\]
Since $\partial_sN=-E$, we set, for $s\ne0$,
\begin{equation}\label{eq:s78-d-quotient}
d_{\varepsilon,a}=\frac{N}{\partial_s f_{\varepsilon,a}}
=\frac{\mathcal K_\varepsilon I_\varepsilon
        +W_\varepsilon\eta_\varepsilon'}
       {a\sqrt q\,\sinh(\sqrt q\,s)
        -W_\varepsilon\eta_\varepsilon'}.
\end{equation}

\textbf{\textit{Smoothness and support.}}
We next show that $d_{\varepsilon,a}$ extends smoothly to $s=0$.
Put $\tau=s/\varepsilon$ and
\[
A(\tau)=\frac1\tau\int_0^\tau\eta(t)\,dt,\qquad
L(\tau)=\frac{\eta'(\tau)}\tau,\qquad
S(v)=\frac{\sinh v}{v},
\]
with smooth values at zero. Indeed, $\eta$ is even,
$\eta=1$ near zero, and $S(v)=\int_0^1\cosh(tv)\,dt\ge1$.
Dividing numerator and denominator by $s$ yields
\begin{align}
\frac Ns&=\mathcal K_\varepsilon A(\tau)
       +\frac{W_\varepsilon}{\varepsilon^2}L(\tau),\notag\\
\frac{\partial_s f_{\varepsilon,a}}s
 &=aqS(\sqrt q\,s)
       -\frac{W_\varepsilon}{\varepsilon^2}L(\tau)
 \ge aq-C\varepsilon^2.
\label{eq:s78-divided}
\end{align}
Here $\|W_\varepsilon\|_\infty\le C\varepsilon^4$.
Since $\varepsilon^2/a\le\delta_*$, choosing
$\delta_*\le q/(2C)$ gives
\[
\frac{\partial_s f_{\varepsilon,a}}s
\ge a(q-C\delta_*)\ge\frac{aq}{2}>0.
\]
Since $(\partial_s f_{\varepsilon,a})/s$ is strictly positive,
the quotient defining $d_{\varepsilon,a}$ extends smoothly to $s=0$.
Consequently, $d_{\varepsilon,a}\partial_s f_{\varepsilon,a}=N$
everywhere. Hence
\[
\operatorname{div}(A_{\varepsilon,a}\nabla f_{\varepsilon,a})
-qf_{\varepsilon,a}=E+\partial_sN=0.
\]
Since $A(0)=1$ and $L(0)=0$, we also obtain
\begin{equation}\label{eq:s78-center-value}
d_{\varepsilon,a}(X',0)=\frac{\mathcal K_\varepsilon(X')}{qa}.
\end{equation}
The zero integral in \eqref{eq:s78-eta} implies
$I_\varepsilon=0$ for $|s|\ge2\varepsilon$; hence $N$ and
$d_{\varepsilon,a}$ have the stated support. When $|X'|<\varepsilon$ and $|s|<\varepsilon$, we have
$W_\varepsilon=P$, $\mathcal K_\varepsilon=0$, and
$\eta_\varepsilon'=0$, so $d_{\varepsilon,a}=0$.

\textbf{\textit{Uniform derivative estimates.}}
For the estimates, set $X'=\varepsilon\xi$, $s=\varepsilon\tau$ and
\[
w(\xi)=\chi(\xi)P(\xi),\qquad k_0=\Delta_\xi w,\qquad
k_1=\frac1{24}\operatorname{div}_\xi
     \bigl(\xi(\xi\cdot\nabla_\xi w)\bigr)-w.
\]
These are fixed smooth compactly supported functions. By the
quartic homogeneity of $P$,
\begin{equation}\label{eq:s78-K-scale}
W_\varepsilon(\varepsilon\xi)=\varepsilon^4w(\xi),\qquad
\mathcal K_\varepsilon(\varepsilon\xi)
 =\varepsilon^2k_0(\xi)+q\varepsilon^4k_1(\xi).
\end{equation}
With $\delta=\varepsilon^2/a$, the divided formulas give
\begin{equation}\label{eq:s78-profile-scaling}
d_{\varepsilon,a}(\varepsilon\xi,\varepsilon\tau)
=\frac{\varepsilon^2}{a}\,
 \mathcal D_{\varepsilon,\delta}(\xi,\tau),
\end{equation}
where
\begin{equation}\label{eq:s78-profile}
\mathcal D_{\varepsilon,\delta}(\xi,\tau)=
\frac{(k_0(\xi)+q\varepsilon^2k_1(\xi))A(\tau)+w(\xi)L(\tau)}
     {qS(\sqrt q\,\varepsilon\tau)-\delta w(\xi)L(\tau)}.
\end{equation}
Decreasing $\delta_*$ if necessary, the denominator is at least
$q/2$ for $0<\varepsilon\le1$ and $0\le\delta\le\delta_*$.
By \eqref{eq:s78-eta}, $A,L$ vanish for $|\tau|\ge2$,
and $w,k_0,k_1$ vanish for $|\xi|\ge2$.
Thus the normalized profiles have fixed compact support and
uniformly bounded derivatives of every fixed order.
Differentiating \eqref{eq:s78-profile-scaling} in physical
coordinates proves \eqref{eq:s78-d-estimate}.
\end{proof}

\subsection{The global construction}\label{sec:s78-global}

We now apply Proposition~\ref{prop:s78-one-patch} at the centers $z_j$
from Corollary~\ref{cor:s56-infinite-perturbations}, with radii
\begin{equation}\label{eq:s78-epsilon-choice}
\varepsilon_j=a_j^2=e^{-2j^4}.
\end{equation}
The small radii allow disjoint corrections while controlling their
first derivatives.

\begin{proposition}\label{prop:s78-global}
For every sufficiently large \(j_0\), there exist a real symmetric,
uniformly elliptic matrix \(A_3\) and a function \(\phi_4\) such that
\[
 A_3\in C^{1,\gamma}(\overline{B_1})\quad
 \text{for every }0<\gamma<1/2,\qquad
 \phi_4\in C^3(\overline{B_1}).
\]
Moreover,
\begin{equation}\label{eq:s78-global-equations}
 K_4\phi_4=\lambda\phi_4,\qquad K_4\psi_2=\mu\psi_2,\qquad
 K_4=-\rho^{-1}{\rm div}(\rho A_3\nabla)+V,
 \qquad \phi_4|_{\partial B_1}=\psi_2|_{\partial B_1}=0.
\end{equation}
On \(\mathcal C_j=\{r<\varepsilon_j,\ |z-z_j|<\varepsilon_j\}\),
the exact identities are
\begin{equation}\label{eq:s78-final-core}
 A_3=A_*,\quad \rho=1,\quad V=\mu,\quad\psi_2=1,\qquad
 \phi_4=a_j\cosh(\sqrt q(z-z_j)).
\end{equation}
In particular, on each two-dimensional disk $D_j=\{z=z_j,\ r<\varepsilon_j\}$,
\[
\phi_4=a_j,\quad\nabla\phi_4=0,
\qquad\psi_2=1,\quad\nabla\psi_2=0.
\]
\end{proposition}

\begin{proof}
\textbf{\textit{Construction.}} By Proposition~\ref{prop:s56-global-cutoff}, we have obtained a smooth real symmetric
uniformly elliptic matrix \(A_2\), smooth \(\rho>0,V\),
and smooth real Dirichlet eigenfunctions satisfying
\[
 K_3\phi_3=\lambda\phi_3,\qquad
 K_3\psi_2=\mu\psi_2,\qquad
 K_3=-\rho^{-1}{\rm div}(\rho A_2\nabla)+V,
\]
where \(q=\mu-\lambda>0\).
For some fixed \(r_0>0\), the central cylinders
\[
 \mathcal C_j^0=\{r<r_0,\ |z-z_j|<2\ell_j/5\}
\]
are interior to \(B_1\), and on them
\begin{equation}\label{eq:s78-global-input}
 A_2=A_*,\quad \rho=1,\quad V=\mu,\quad\psi_2=1,\qquad
 \phi_3=P+a_j\cosh(\sqrt q(z-z_j)),
\end{equation}
with
\[
 z_j=2^{-j},\qquad \ell_j=2^{-j-4},\qquad a_j=e^{-j^4}.
\]
The choice of radii gives
\[
 \frac{\varepsilon_j}{\ell_j}\to0,\qquad
 \frac{\varepsilon_j}{z_j}\to0,\qquad
 \frac{\varepsilon_j^2}{a_j}=a_j^3\to0.
\]
Increase \(j_0\) so that Proposition~\ref{prop:s78-one-patch} applies
and, for every \(j\ge j_0\),
\[
 2\varepsilon_j<r_0,\qquad
 2\varepsilon_j<2\ell_j/5,\qquad
 2\varepsilon_j<z_j/4.
\]
The closed support cylinders
\[
 S_j=\{r\le2\varepsilon_j,\ |z-z_j|\le2\varepsilon_j\}
\]
are contained in \(\mathcal C_j^0\) and are pairwise disjoint,
because their longitudinal projections lie inside the disjoint
intervals \(I_j=(z_j-\ell_j,z_j+\ell_j)\).
Their only accumulation point is \(0\), and on \(S_j\),
\begin{equation}\label{eq:s78-distance}
 c z_j\le |X|\le C z_j
\end{equation}
for fixed positive constants \(c,C\).

Applying Proposition~\ref{prop:s78-one-patch} with
$\varepsilon=\varepsilon_j$, $a=a_j$, and $s=z-z_j$, define
\[
 \Gamma_j(X)=
 W_{\varepsilon_j}(X')\,
 \eta((z-z_j)/\varepsilon_j),\qquad
 d_j(X)=d_{\varepsilon_j,a_j}(X',z-z_j),
\]
using smooth zero extension.
Set
\begin{equation}\label{eq:s78-global-definitions}
 \Gamma=\sum_{j\ge j_0}\Gamma_j,\qquad
 d=\sum_{j\ge j_0}d_j,\qquad
 \phi_4=\phi_3-\Gamma.
\end{equation}
Define the assembled matrix \(A_3\) entry by entry:
\begin{equation}\label{eq:s78-global-B}
 (A_3)_{33}=(A_2)_{33}+d,\qquad
 (A_3)_{ik}=(A_2)_{ik}\quad\text{if }(i,k)\ne(3,3).
\end{equation}
There is at most one nonzero summand at any point.
At \(0\), all corrections are defined to be zero.
The density, potential, and second eigenfunction are unchanged.

\textbf{\textit{Ellipticity and $C^1$ regularity.}}
By \eqref{eq:s78-d-estimate} and \eqref{eq:s78-epsilon-choice},
\begin{equation}\label{eq:s78-low-bounds}
 ||d_j||_\infty\le Ca_j^3,\qquad
 ||\nabla d_j||_\infty\le Ca_j,\qquad
 ||D^2d_j||_\infty\le C/a_j.
\end{equation}
Here \(D^2d_j\) denotes the array of second coordinate derivatives.
If \(\lambda_0>0\) is an ellipticity lower bound for
\(A_2\), disjointness gives
\[
 ||d||_\infty
 \le C\sup_{j\ge j_0}a_j^3.
\]
Choose \(j_0\) so that this is below \(\lambda_0/2\).
For every \(\xi\in\R^3\),
\[
 \xi^T A_3(X)\xi
 =\xi^TA_2(X)\xi+d(X)\xi_3^2
 \ge\frac{\lambda_0}{2}|\xi|^2.
\]
The upper bound and symmetry are also preserved. Away from \(0\), the sums are locally finite and smooth.
For \(X\in S_j\), \eqref{eq:s78-distance} gives
\[
 \frac{|d(X)|}{|X|}
 \le C\frac{a_j^3}{z_j}\longrightarrow0.
\]
Outside the support cylinders, the correction is zero.
This proves differentiability at \(0\), with
\(\nabla d(0)=0\).
The gradient estimates in \eqref{eq:s78-low-bounds} give
\(\nabla d(X)\to0\) as \(X\to0\).
Thus \(d\in C^1\).

\textbf{\textit{The global H\"older estimate.}}
Use the smooth zero extension of each \(d_j\) to \(\R^3\).
The mean value theorem and \eqref{eq:s78-low-bounds} imply
\[
 |\nabla d_j(X)-\nabla d_j(Y)|
 \le\min\{Ca_j,\;Ca_j^{-1}|X-Y|\}.
\]
For positive \(a,b\) and \(0<\gamma<1\),
\(\min\{a,b\}\le a^{1-\gamma}b^\gamma\).
Hence
\begin{equation}\label{eq:s78-holder}
 [\nabla d_j]_{C^{0,\gamma}(\R^3)}
 \le C_\gamma a_j^{1-2\gamma}.
\end{equation}
For \(0<\gamma<1/2\), these constants are uniformly bounded.
If \(X,Y\) lie in one support cylinder, apply the estimate to
that correction. If they lie in different cylinders \(S_j,S_k\),
then \(\nabla d_j(Y)=\nabla d_k(X)=0\), and therefore
\begin{align*}
 |\nabla d(X)-\nabla d(Y)|
 &\le|\nabla d_j(X)-\nabla d_j(Y)|\\
 &\quad+|\nabla d_k(X)-\nabla d_k(Y)|\\
 &\le 2\sup_\ell[\nabla d_\ell]_{C^{0,\gamma}}|X-Y|^\gamma.
\end{align*}
If one point lies outside all cylinders, the same argument uses
one correction. It also applies when one point is \(0\).
This proves the asserted global \(C^{1,\gamma}\) regularity.

\textbf{\textit{Regularity of the modified eigenfunction.}}
Recall from the proof of Proposition~\ref{prop:s78-one-patch}
that $w(\xi)=\chi(\xi)P(\xi)$. Since $w$ and $\eta$ are fixed
smooth compactly supported functions, the scaling identity
\eqref{eq:s78-K-scale} gives
\[
 \Gamma_j(X)
 =\varepsilon_j^4
 w(X'/\varepsilon_j)\eta((z-z_j)/\varepsilon_j).
\]
Thus, for every multi-index \(\alpha\),
\begin{equation}\label{eq:s78-Gamma-estimate}
 ||\partial^\alpha\Gamma_j||_\infty
 \le C_\alpha\varepsilon_j^{4-|\alpha|}.
\end{equation}
For \(|\alpha|\le3\), these quantities tend to zero.
For $|\alpha|\le3$, extend $\partial^\alpha\Gamma$,
initially defined away from the origin, by zero at $0$.
The estimates show that these extensions are continuous.
To verify actual differentiability inductively, for
an order \(m-1\) derivative, \(1\le m\le3\), use
\[
 \frac{|\partial^\alpha\Gamma(X)|}{|X|}
 \le C\frac{\varepsilon_j^{5-m}}{z_j}\longrightarrow0,
 \qquad X\in S_j,\quad |\alpha|=m-1.
\]
The limit follows from
\(\varepsilon_j=e^{-2j^4}\), \(z_j=2^{-j}\).
Outside the cylinders the numerator is zero.
This proves \(\Gamma\in C^3\), with all its
derivatives through order three zero at \(0\).
Consequently \(\phi_4\in C^3(\overline{B_1})\).

\textbf{\textit{Verification of the equations.}}
On $S_j$, \eqref{eq:s78-global-input} gives $A_2=A_*$,
$\rho=1$, $V=\mu$, and
$\phi_3=P+a_j\cosh(\sqrt q(z-z_j))$. Moreover,
\[
\phi_4(X)=f_{\varepsilon_j,a_j}(X',z-z_j),\qquad
A_3(X)=A_{\varepsilon_j,a_j}(X',z-z_j).
\]
Hence Proposition~\ref{prop:s78-one-patch} gives
\({\rm div}(A_3\nabla\phi_4)=q\phi_4\), equivalently
\(K_4\phi_4=\lambda\phi_4\).
Away from $0$ and the support cylinders, $A_3=A_2$ and
$\phi_4=\phi_3$, so $K_4\phi_4=\lambda\phi_4$.
The matrix changes occur where \(\psi_2=1\), and hence do not
alter the second eigenfunction equation for \(\psi_2\).
By continuity, both eigenvalue equations also hold at the origin.
All changes have compact interior support, so the boundary traces
remain zero. This proves \eqref{eq:s78-global-equations}.

Finally, on $\mathcal C_j$, Proposition~\ref{prop:s78-one-patch}
gives $\phi_4=a_j\cosh(\sqrt q(z-z_j))$ and $d_j=0$.
All other corrections vanish there by disjointness.
Together with \eqref{eq:s78-global-input}, this proves
\eqref{eq:s78-final-core}; the disk identities follow by setting $z=z_j$.
\end{proof}

\subsection{Failure of higher regularity}

We examine the higher-order derivatives of the local corrections as
their supports approach the origin.

\begin{proposition}\label{prop:s78-regularity-loss}
The assembled coefficient \(A_3\) is not \(C^2\) in any neighborhood
of \(0\), and \(\phi_4\) is not \(C^4\) in any neighborhood of \(0\).
\end{proposition}

\begin{proof}
\textbf{\textit{Failure of $C^2$ regularity of $A_3$.}}
Let \(w,k_1\) be the fixed functions in \eqref{eq:s78-K-scale}.
By \eqref{eq:s78-center-value}, for \(i,k\in\{1,2\}\),
\begin{equation}\label{eq:s78-second-derivative}
\partial_{X_i}\partial_{X_k}d_j(\varepsilon_j\xi,z_j)
=\frac{\partial_{\xi_i}\partial_{\xi_k}\Delta_\xi w(\xi)}{qa_j}
+\frac{\varepsilon_j^2}{a_j}\partial_{\xi_i}\partial_{\xi_k}k_1(\xi).
\end{equation}
Some entry of the Hessian of \(\Delta_\xi w\) is nonzero at a
fixed point \(\xi_*\).
Indeed, if all its entries vanished identically, \(\Delta_\xi w\)
would be affine. It is compactly supported, so it would be zero.
Integration by parts would then give
\[
 \int_{\R^2}|\nabla w|^2\,d\xi
 =-\int_{\R^2}w\Delta_\xi w\,d\xi=0,
\]
forcing \(w=0\), contrary to \(w=P\) on the unit disk.
Fix such an entry and point, and write its nonzero value as \(c_*\).
Since \(\varepsilon_j\to0\), formula
\eqref{eq:s78-second-derivative} implies, for large \(j\),
\[
 \left|
 \partial_{X_i}\partial_{X_k}d_j(\varepsilon_j\xi_*,z_j)
 \right|\ge\frac{|c_*|}{2qa_j}\longrightarrow\infty.
\]
These points tend to \(0\). Since $|\xi_*|<2$, they belong
to $S_j\subset\mathcal C_j^0$, where \((A_2)_{33}=1\)
on an open neighborhood.
Other corrections vanish there, so the same diverging derivatives
are derivatives of \((A_3)_{33}\).
This is incompatible with continuity of second derivatives near \(0\).

\textbf{\textit{Failure of $C^4$ regularity of $\phi_4$.}}
For the eigenfunction, set \(p_j=(0,0,z_j)\).
By \eqref{eq:s78-final-core},
\[
 \partial_x^4\phi_4(p_j)=0.
\]
For each fixed \(|X'|<r_0\), continuity of \(\phi_3\)
and \eqref{eq:s78-global-input} imply
\[
 \phi_3(X',0)
 =\lim_{j\to\infty}\left(P(X')+a_j\right)=P(X').
\]
Since $\operatorname{supp}\Gamma_j\subset\{z>0\}$ for every $j$,
we have $\Gamma(X',0)=0$.
Hence \(\phi_4(x,y,0)=P(x,y)\) for \(r<r_0\).
If \(\phi_4\) were \(C^4\) near \(0\), this identity would imply
\(\partial_x^4\phi_4(0)=24\).
Continuity along \(p_j\to0\) would instead give zero, a contradiction.
\end{proof}
% --------------------------------------------------------------------
% Removing the potential and density
\section{Completion of the counterexample construction}\label{removing-the-potential-and-density}

Proposition~\ref{prop:s78-global} provides two Dirichlet eigenfunctions
for an operator with density \(\rho\) and potential \(V\).
We remove these terms by a positive conjugation and a change of
coordinates to obtain the heat flow in Theorem~\ref{thm:counter}.
The following theorem restates Theorem~\ref{thm:counter} in a more explicit form.

\begin{theorem}\label{thm:counter-construction}
For every $0<\gamma<1/2$, there exist a real symmetric,
time-independent, uniformly elliptic matrix
$A\in C^{1,\gamma}(\overline{B_1})$, two nonzero Dirichlet
eigenfunctions $\Phi,\Psi$ of $-\operatorname{div}(A\nabla)$
with distinct eigenvalues $\lambda'<\mu'$, and a diffeomorphism
$F:\overline{B_1}\to\overline{B_1}$ such that
\[
u(Y,t)=e^{-\lambda't}\Phi(Y)-e^{-\mu't}\Psi(Y)
\]
is a nonzero solution of
\[
\partial_tu=\operatorname{div}(A\nabla u)
\quad\text{in }B_1,\qquad u|_{\partial B_1}=0.
\]
There are pairwise disjoint two-dimensional disks $D_j\subset B_1$
and times $t_j\to\infty$ for which
\[
u(Y,t_j)=0,\qquad \nabla_Yu(Y,t_j)=0
\quad\text{for every }Y\in F(D_j).
\]
In particular, $\dim_{\mathrm H}S_{t_j}(u)\ge2$.
\end{theorem}

\begin{proof}
\textbf{\textit{Removing the potential.}}
The operator from Proposition~\ref{prop:s78-global} is
\[
K_4=-\rho^{-1}\operatorname{div}(\rho A_3\nabla)+V,
\]
where \(\rho>0,V\) are smooth and \(A_3\) is symmetric and
uniformly elliptic. Fix \(0<\gamma<1/2\). By
Proposition~\ref{prop:s78-global}, \(A_3\in C^{1,\gamma}(\overline{B_1})\).
Choose a common spectral shift \(C\) such that \(V+C>0\), and solve
for \(p\):
\begin{equation*}
(K_4+C)p=0\quad\text{in }B_1,
\qquad p|_{\partial B_1}=1.
\end{equation*}
Elliptic Schauder estimates give \(p\in C^{2,\gamma}(\overline{B_1})\),
while the maximum principle gives \(p>0\) on \(\overline{B_1}\).

Recall the auxiliary solution \(U\) from the start of Part II, now with \(\phi=\phi_4\) and \(\psi=\psi_2\). Set \(\lambda'=\lambda+C\) and \(\mu'=\mu+C\), and define the shifted
two-mode solution
\begin{equation*}
U_C(X,t)=e^{-\lambda' t}\phi_4(X)-e^{-\mu' t}\psi_2(X)
      =e^{-Ct}U(X,t).
\end{equation*}
Then $\partial_tU_C+(K_4+C)U_C=0$. Since
$U_C=e^{-Ct}U$, the function $U_C$ and its spatial gradient
vanish on $D_j$ at the times
\[
t_j=q^{-1}\log(1/a_j)=j^4/q\longrightarrow\infty.
\]
Define the conjugated eigenfunction pair by \(\phi_5=\phi_4/p\) and
\(\psi_3=\psi_2/p\).
Set \(v=e^{-\lambda' t}\phi_5-e^{-\mu' t}\psi_3=U_C/p\) and
\(\sigma=\rho p^2\).
Direct calculation of the conjugation by a positive function yields
\begin{equation}
v_t=\sigma^{-1}\operatorname{div}(\sigma A_3\nabla v),
\qquad v|_{\partial B_1}=0.
\label{eq:density-heat-equation}
\end{equation}
\textbf{\textit{Removing the density.}}
Since \(\sigma=\rho p^2\) is positive and belongs to
\(C^{2,\gamma}(\overline{B_1})\), set
\begin{equation*}
c_0=|B_1|^{-1}\int_{B_1}\sigma\,dX.
\end{equation*}
The Dacorogna--Moser prescribed-Jacobian theorem \cite{dm1990}
gives a diffeomorphism \(F:\overline{B_1}\to\overline{B_1}\), with
\(F,F^{-1}\in C^{3,\gamma}(\overline{B_1})\), satisfying
\begin{equation*}
\det DF=\sigma/c_0,\qquad F|_{\partial B_1}=\mathrm{id}.
\end{equation*}
In the coordinates \(Y=F(X)\), set \(\widetilde v(Y,t)=v(F^{-1}(Y),t)\).
Let $\eta$ be a smooth compactly supported test function and set
$\widetilde\eta=\eta\circ F^{-1}$. Multiplying
\eqref{eq:density-heat-equation} by $\sigma\eta$, integrating by parts,
and changing variables from $X$ to $Y=F(X)$, we obtain
\begin{equation*}
\int_{B_1}\sigma v_t\eta\,dX
=c_0\int_{B_1}\widetilde v_t\widetilde\eta\,dY,
\end{equation*}
\begin{equation*}
\int_{B_1}\sigma A_3\nabla v\cdot\nabla\eta\,dX
=c_0\int_{B_1}
(DF\,A_3\,DF^T)\nabla\widetilde v\cdot\nabla\widetilde\eta\,dY.
\end{equation*}
Canceling the common constant \(c_0\) gives
\begin{equation*}
\partial_t\widetilde v
=\operatorname{div}_Y(A(Y)\nabla_Y\widetilde v),
\end{equation*}
\begin{equation*}
A(Y)=
\big(DF(X)A_3(X)DF(X)^T\big)\big|_{X=F^{-1}(Y)}.
\end{equation*}
The matrix \(A\) is time-independent, symmetric, and uniformly elliptic.
Since $A_3\in C^{1,\gamma}(\overline{B_1})$ and
$F,F^{-1}\in C^{3,\gamma}(\overline{B_1})$, its defining formula gives
$A\in C^{1,\gamma}(\overline{B_1})$.

Define the final Dirichlet eigenfunctions by
\(\Phi=\phi_5\circ F^{-1}\) and \(\Psi=\psi_3\circ F^{-1}\).
The boundary of the ball is fixed, and the final heat flow
\(u=\widetilde v\) satisfies the Dirichlet condition. Explicitly,
\begin{equation}\label{eq:counterexample-final-flow}
u(Y,t)=e^{-\lambda' t}\Phi(Y)-e^{-\mu' t}\Psi(Y),\qquad
u(F(X),t)=\frac{e^{-\lambda' t}\phi_4(X)-e^{-\mu' t}\psi_2(X)}{p(X)}.
\end{equation}

\textbf{\textit{Critical disks.}}
Multiplication by a positive function preserves the critical zero set:
at points where \(U_C=0\), we have \(\nabla(U_C/p)=\nabla U_C/p\). A
diffeomorphism also preserves vanishing of the gradient and Hausdorff
dimension. Consequently,
\begin{equation*}
F(D_j)\subset S_{t_j}(u),
\qquad t_j=q^{-1}\log(1/a_j)=j^4/q\longrightarrow\infty,
\end{equation*}
and
\begin{equation*}
\dim_H S_{t_j}(u)\ge2.
\end{equation*}
The solution \(u\) is nonzero, since it arises from two
nonzero modes with distinct eigenvalues.

This completes the proof of Theorem~\ref{thm:counter}.
\end{proof}

% --------------------------------------------------------------------
% Regularity thresholds and open questions

\end{document}